\documentclass[11pt]{article}

\usepackage[utf8]{inputenc}

\usepackage{mathtools}
\mathtoolsset{centercolon}
\usepackage{microtype}
\usepackage[usenames,dvipsnames,svgnames,table]{xcolor}
\usepackage{amsfonts,amsmath, amssymb,amsthm,amscd,bm}
\usepackage[height=9in,width=6.5in]{geometry}
\usepackage{verbatim}
\usepackage{tikz}
\usetikzlibrary{decorations.markings}
\tikzset{->-/.style={decoration={
			markings,
			mark=at position .6 with {\arrow{>}}},postaction={decorate}}}
\usepackage{tkz-euclide}
\usetikzlibrary{arrows.meta}

\usepackage[margin=10pt,font=small,labelfont=bf]{caption}
\usepackage{latexsym}
\usepackage[pdfauthor={Tang},bookmarksnumbered,hyperfigures,colorlinks=true,citecolor=BrickRed,linkcolor=BrickRed,urlcolor=BrickRed,pdfstartview=FitH]{hyperref}
\usepackage{mathrsfs}
\usepackage{graphicx}
\usepackage{enumitem}

\usepackage[inline,nolabel]{showlabels}

\newtheorem{theorem}{Theorem}[section]

\newtheorem{lemma}[theorem]{Lemma}
\newtheorem{proposition}[theorem]{Proposition}
\newtheorem{question}[theorem]{Question}
\newtheorem{conj}[theorem]{Conjecture}

\newtheorem{definition}[theorem]{Definition}
\newtheorem{remark}[theorem]{Remark}
\newtheorem{example}[theorem]{Example}
\numberwithin{equation}{section}

\newcommand{\be}{\begin{equation}}
	\newcommand{\ee}{\end{equation}}
\newcommand\ba{\begin{align}}
	\newcommand\ea{\end{align}}

\newcommand{\notion}[1]{{\bf  \textit{#1}}}

\begin {document}
\title{A note on some critical thresholds of Bernoulli percolation}
\author{
	Pengfei Tang\thanks{Department of Mathematical Sciences, Tel Aviv
		University.  Partially supported by the National
		Science Foundation under grant DMS-1612363, DMS-1954086 and by ERC starting grant 676970 RANDGEOM.
		Email: \protect\url{pengfeitang@mail.tau.ac.il}.}
}
\date{\today}
\maketitle

% MR: 60J10  	Markov chains (discrete-time Markov processes on discrete state spaces)
% 05C81  	Random walks on graphs
% 05C05  	Trees

% keywords: trees; cover time; BEST theorem

\begin{abstract}
Consider Bernoulli bond percolation on a locally finite, connected graph $G$ and let $p_{\mathrm{cut}}$ be the threshold corresponding to a ``first-moment method" lower bound.   Kahn (\textit{Electron.\ Comm.\ Probab.\ Volume 8, 184-187.} (2003)) constructed a counter-example to Lyons' conjecture of $p_{\mathrm{cut}}=p_c$ and proposed a modification. Here we give a positive answer to Kahn's modified question. The key observation is that in Kahn's modification, the new expectation quantity also appears in the differential inequality of one-arm events. This links the question to a lemma of Duminil-Copin and Tassion (\textit{Comm. Math. Phys. Volume 343, 725-745.} (2016)). We also study some applications for Bernoulli percolation on periodic trees. 
	
%For a locally finite, connected graph $G=(V,E)$, let $p_c=p_c(G)$ denote the critical probability for Bernoulli bond percolation on $G$. Consider Bernoulli$(p)$ bond percolation on $G$. For $x\in V$ and a finite set $S$ with $x\in S$, let $\partial_E S$ denote the edge boundary of $S$ and let $\varphi_p(x,S)$ denote the expected number of open edges in $\partial_E S$ that are connected to $x$ via an open path in $S$. In the new proof of sharpness of phase transitions for Bernoulli percolation and the Ising model, Duminil-Copin and Tassion (\textit{Comm. Math. Phys. Volume 343, 725-745.} (2016)) showed that for transitive graphs,
%\[
%p_c=\sup\{p\geq0\colon \inf\{\varphi_p(x,S)\colon x\in S, \, S \textnormal{ is finite}\} <1 \}.
%\] 
%
%In this note we point out that for every locally finite, connected, infinite graph, 
%\[
%p_c=\sup\{p\geq0\colon \inf\{\varphi_p(x,S)\colon x\in S, \, S \textnormal{ is finite}\} =0 \}
%\]
%and this resolves a question of Kahn (\textit{Electron. Comm. Probab. Volume 8, 184-187.} (2003)).
\end{abstract}

\section{Introduction}\label{sec:intro}

Let $G=(V,E)$ be a locally finite (i.e., each vertex has finite degree),   connected, infinite graph. For $p\in[0,1]$, \notion{Bernoulli$(p)$ bond percolation} studies the random subgraph $\omega$ of $G$ formed by keeping each edge with probability $p$ and removing otherwise, independently of each other. The edges kept in $\omega$ are called \notion{open edges} and the edges removed are called \notion{closed edges}.   The connected components are called (open) clusters.  For background on Bernoulli percolation, see Chapter 7 of \cite{LP2016} or \cite{Grimmett1999}. For $p\in[0,1]$, let $\mathbb{P}_p$  denote the law of Bernoulli$(p)$ bond percolation and $\mathbb{E}_p$ the corresponding expectation. 

% One can also study the induced subgraph by keeping each vertex with probability $p$ and removing  otherwise, independently of each other. We call the random subgraph obtained in this way   Bernoulli$(p)$ site percolation and denote its law by $\mathbb{P}_p^{\mathrm{site}}$. 

Let $C(x)$ denote the open cluster of $x$ in Bernoulli percolation. Let $|C(x)|_V,|C(x)|_E$ denote the number of vertices and edges in the cluster $C(x)$ respectively. Let $A\longleftrightarrow B$ denote the event that there is an open path connecting some vertex $x\in A$ and $y\in B$. Let $x\longleftrightarrow \infty$  denote the event that the diameter of $C(x)$ is infinite. The \notion{critical probability} $p_{\mathrm{c}}$ is defined as 
\[
p_{\mathrm{c}}=p_{\mathrm{c}}(G):=\sup\{p\geq0\colon \mathbb{P}_p(x\longleftrightarrow \infty)=0 \}.
\]

Since for a locally finite graph $G$, the three events $x\longleftrightarrow\infty$, $|C(x)|_V=\infty$ and $|C(x)|_E=\infty$ are actually the same event, one can also define 
\[
p_{\mathrm{c}}=p_{\mathrm{c}}(G):=\sup\{p\geq0\colon \mathbb{P}_p(|C(x)|_V=\infty)=0 \}
\]
or 
\[
p_{\mathrm{c}}=p_{\mathrm{c}}(G):=\sup\{p\geq0\colon \mathbb{P}_p(|C(x)|_E=\infty)=0 \}.
\]

Let $x$ be a vertex in $G$. We say that $\Pi_E$ is a \notion{edge cutset} separating $x$ from infinity, if $\Pi_E$ is a set of edges such that the connected component of $x$ in $G\backslash \Pi_E$ is finite. Similarly one can define \notion{vertex cutsets}.

\begin{definition}\label{def: p_cut,E and p_cut,V}
	Suppose $G$ is a locally finite, connected, infinite graph. Define 
	\begin{equation*}\label{eq: def of p_cut,E in bond percolation}
		p_{\mathrm{cut,E}}= p_{\mathrm{cut,E}}(G):=\sup\{p\geq0\colon \inf_{\Pi_E}\mathbb{E}_p[|C(x)\cap \Pi_E|]=0 \},
	\end{equation*}
	where the infimum is taken over all \textbf{edge cutsets} $\Pi_E$ separating $x$ from infinity and $C(x)\cap \Pi_E$ denotes the intersection of the edge set of $C(x)$ with $\Pi_E$.

	Define 	
	\begin{equation*}\label{eq: def of p_cut,V in bond percolation}
		p_{\mathrm{cut,V}}= p_{\mathrm{cut,V}}(G):=\sup\{p\geq0\colon \inf_{\Pi_V}\mathbb{E}_p[|C(x)\cap \Pi_V|]=0 \},
	\end{equation*}
	where the infimum is taken over all \textbf{vertex cutsets}  $\Pi_V$ separating $x$ from infinity and $C(x)\cap \Pi_V$ denotes the intersection of the vertex set of $C(x)$ with $\Pi_V$.  	
\end{definition}

For any edge (or vertex) cutset $\Pi$ separating $x$ from infinity, if the event $\{x\longleftrightarrow\infty\}$ occurs, then $C(x)\cap\Pi$ is nonempty. Hence 
\be\label{eq: p_cut leq p_c's derivation}
\mathbb{P}_p(x\longleftrightarrow \infty)\leq \mathbb{P}_p(|C(x)\cap \Pi|\geq 1)\leq \mathbb{E}_p[|C(x)\cap \Pi|].
\ee
Thus one has that 
\be \label{eq: p_cut leq p_c}
p_{\mathrm{cut,E}}\leq p_{\mathrm{c}}\,\,\,\textnormal{ and }\,\,\,	p_{\mathrm{cut,V}}\leq p_{\mathrm{c}}.
\ee

Historically another critical value $p_T$ is also of great interest (coincide with the notation $p_{\mathrm{T,V}}$ below). 
\begin{definition}
	Suppose $G$ is a locally finite, connected, infinite graph. 
	Define 
	\begin{equation*}\label{eq: def of p_T,V in bond percolation}
		p_{\mathrm{T,V}}=p_{\mathrm{T,V}}(G):=\sup\{p\geq0\colon \mathbb{E}_p[|C(x)|_V]<\infty \}
	\end{equation*}
	and 
	\begin{equation*}\label{eq: def of p_T,E in bond percolation}
		p_{\mathrm{T,E}}=p_{\mathrm{T,E}}(G):=\sup\{p\geq0\colon \mathbb{E}_p[|C(x)|_E]<\infty \}.
	\end{equation*}
\end{definition}

If $p<p_{\mathrm{T,V}}$, then $\sum_{n=1}^{\infty}\mathbb{E}_p[|C(x)\cap \Pi_n|]\leq  \mathbb{E}_p[|C(x)|_V]<\infty$, where $\Pi_n:=\{y\colon d_G(y,x)=n \}$ is the cutset consisting of vertices at graph distance $n$ to $x$. Hence $p<p_{\mathrm{T,V}}$ implies that $p\leq p_{\mathrm{cut,V}}$. Thus 
\be\label{eq: p_T,V leq p_cut,V}
p_{\mathrm{T,V}}\leq p_{\mathrm{cut,V}}.
\ee
Similarly one has that 
\be\label{eq: p_T,E leq p_cut,E}
p_{\mathrm{T,E}}\leq p_{\mathrm{cut,E}}.
\ee

It is easy to see that these critical values $p_{\mathrm{c}},p_{\mathrm{cut,E}},p_{\mathrm{cut,V}},p_{\mathrm{T,E}},p_{\mathrm{T,V}}$ do not depend on the choice of $x$ by Harris' inequality \cite[Section 5.8]{LP2016}. 

By \eqref{eq: p_cut leq p_c}, \eqref{eq: p_T,V leq p_cut,V} and \eqref{eq: p_T,E leq p_cut,E} we now have 
\[
p_{\mathrm{T,E}}\leq p_{\mathrm{cut,E}}\leq p_{\mathrm{c}}.
\]
and
\[
p_{\mathrm{T,V}}\leq p_{\mathrm{cut,V}}\leq p_{\mathrm{c}}.
\]

Lyons showed that $p_{\mathrm{c}}=p_{\mathrm{cut,V}}$ holds for trees \cite{Lyons1990} and tree-like graphs \cite{Lyons1989} and pointed out $p_{\mathrm{cut,V}}=p_{\mathrm{c}}$ for transitive graphs in \cite{Lyons1990} since $p_{\mathrm{T,V}}=p_{\mathrm{c}}$ for such graphs \cite{AB1987,Men1986}; and these results for $p_{\mathrm{cut,V}}$ applied equally to $p_{\mathrm{cut,E}}$ on these graphs. In view of these examples Lyons conjectured that  $p_{\mathrm{c}}=p_{\mathrm{cut,V}}$ for general graphs (lines 11--12 on page 955 of \cite{Lyons1990}). 

Later Kahn \cite{Kahn2003} constructed  a family of counterexamples to Lyons' conjecture.  Kahn's examples exhibited a sequences of vertex cutsets $\Pi_n$ such that  the quantity $|C(x)\cap \Pi_n|=\sum_{v\in \Pi_n}\mathbf{1}_{(x\longleftrightarrow v)}$ is usually zero but has a large expectation for some $p<p_{\mathrm{c}}(G)$. That was achieved by large correlation among the events $\{x\longleftrightarrow v\}$ for $v\in \Pi_n$, i.e., conditioned on the event that $v$ is connected to $x$ via an open path, with high probability a lot of other vertices in $\Pi_n$ are also connected to $x$ via $v$. In light of this Kahn  proposed the following modification of Lyons' conjecture: 
\begin{question}\label{ques: Kahn}
	Does $p_{\mathrm{c}}(G)= p'_{\mathrm{cut,V}}(G)$ hold for every locally finite, connected, infinite graph $G$?
\end{question}
Here  the notation $p'_{\mathrm{cut,V}}=p'_{\mathrm{cut,V}}(G)$ from  \cite{Kahn2003} (there it was denoted by $p'_{\mathrm{cut}}$) is defined as follows.

\begin{definition}\label{def: def of p'_cut,E and p'_cut V}
	Suppose $G$ is a locally finite, connected, infinite graph. For  $x\in V$, let $\Pi_V$ be a vertex cutset which separates $x$ from infinity. For each $v\in\Pi_V$,  let  $A(x,v,\Pi_V)$ denote the event that $x$ is connected to $v$ via an open path without using vertices in $\Pi_V\backslash\{v\}$.
	Define 	
	\begin{equation*}\label{eq: def of p'_cut,V in bond percolation}
		p'_{\mathrm{cut,V}}= p'_{\mathrm{cut,V}}(G):=\sup\Big\{p\geq0\colon \inf_{\Pi_V}\sum_{v\in \Pi_V}\mathbb{P}_p[A(x,v,\Pi_V)]=0 \Big\},
	\end{equation*}
	where the infimum is taken over all \textbf{vertex cutsets}  $\Pi_V$ separating $x$ from infinity. 
	
	Similarly for an edge cutset $\Pi_E$ separating $x$ from infinity and $e\in\Pi_E$, let $A(x,e,\Pi_E)$ denote the event that $x$ is connected to $e$ via an open path without using edges in $\Pi_E\backslash \{e\}$ (Here we assume $e$ itself is also open on $A(x,e,\Pi_E)$.) Define 
	\begin{equation*}\label{eq: def of p'_cut,E in bond percolation}
		p'_{\mathrm{cut,E}}= p'_{\mathrm{cut,E}}(G):=\sup\Big\{p\geq0\colon \inf_{\Pi_E}\sum_{e\in \Pi_E}\mathbb{P}_p[A(x,e,\Pi_E)]=0 \Big\},
	\end{equation*}
	where the infimum is taken over all \textbf{edge cutsets} $\Pi_E$ separating $x$ from infinity. 
	
\end{definition}

Similarly one can ask (\cite[Question 5.16]{LP2016}):
\begin{question}\label{ques: Kahn's question for edge cutset}
	Does $p_{\mathrm{c}}(G)= p'_{\mathrm{cut,E}}(G)$ hold for every locally finite, connected, infinite graph $G$?
\end{question}

Our main result is the following  affirmative answer to Question \ref{ques: Kahn} and \ref{ques: Kahn's question for edge cutset} for Bernoulli bond percolation. 
\begin{theorem}\label{thm: p_c=p_cut'}	
	For Bernoulli bond percolation on every locally finite, connected, infinite graph $G$, one has that 
	\[
	p'_{\mathrm{cut,E}}=	p'_{\mathrm{cut,V}}=p_{\mathrm{c}}.
	\]
\end{theorem}

The same result holds for Bernoulli site percolation on a locally finite, connected, infinite graph with \textbf{bounded degree} if one defines  $ p'_{\mathrm{cut,E}},	p'_{\mathrm{cut,V}}$ accordingly using Bernoulli  site percolation; see Remark \ref{rem: site percolation} and Conjecture \ref{conj: site percolation}  for more discussions.

\section{Some relations of the critical thresholds}\label{sec: simple relations}

For any edge cutset $\Pi$ separating $x$ from infinity, if the event $\{x\longleftrightarrow\infty\}$ occurs, then there is at least one edge $e$ such that the event $A(x,e,\Pi)$ occurs. Hence by union bounds,
\be\label{eq: p'_cut leq p_c's derivation}
\mathbb{P}_p(x\longleftrightarrow \infty)\leq\sum_{e\in \Pi}\mathbb{P}_p[A(x,e,\Pi)]
\ee
Thus one has that 
\be \label{eq: p'_cut,E leq p_c}
p'_{\mathrm{cut,E}}\leq p_{\mathrm{c}}.
\ee
Similarly one has that 
\be \label{eq: p'_cut,V leq p_c}
p'_{\mathrm{cut,V}}\leq p_{\mathrm{c}}.
\ee
Also obviously for any edge cutset $\Pi$ one has that $\sum_{e\in \Pi}\mathbb{P}_p[A(x,e,\Pi)]\leq\sum_{e\in \Pi}\mathbb{P}_p[e\in C(x)]= \mathbb{E}_p[|C(x)\cap \Pi|]$. Hence one has that 
\be\label{eq: p_cut,E less than p'_cut,E}
p_{\mathrm{cut,E}}\leq p'_{\mathrm{cut,E}}.
\ee
Similarly one has that 
\be\label{eq: p_cut,V less than p'_cut,V}
p_{\mathrm{cut,V}}\leq p'_{\mathrm{cut,V}}.
\ee
By \eqref{eq: p_T,V leq p_cut,V}, \eqref{eq: p_T,E leq p_cut,E}, \eqref{eq: p'_cut,E leq p_c},\eqref{eq: p'_cut,V leq p_c}, \eqref{eq: p_cut,E less than p'_cut,E} and \eqref{eq: p_cut,V less than p'_cut,V} one has that 
\be\label{eq: sum of the relations in vertex cutset}
p_{\mathrm{T,V}}\leq p_{\mathrm{cut,V}}\leq p'_{\mathrm{cut,V}} \leq p_{\mathrm{c}}
\ee
and 
\be\label{eq: sum of the relations in edge cutset}
p_{\mathrm{T,E}}\leq p_{\mathrm{cut,E}}\leq p'_{\mathrm{cut,E}} \leq p_{\mathrm{c}}.
\ee

We also have the following relations. 
\begin{lemma}\label{lem: p_cut,e less than p_cut,v}
	Suppose $G$ is a locally finite, connected, infinite graph. 
	Then 
	\be\label{eq: relation of p_cut for edge cutset and vertex cutset}
	p_{\mathrm{cut,E}}\leq  p_{\mathrm{cut,V}}
	\ee
	If moreover $G$ has bounded degree, then the equality holds in \eqref{eq: relation of p_cut for edge cutset and vertex cutset}.
\end{lemma}

\begin{lemma}\label{lem: p_T,E less than P_T,V}
	Suppose $G$ is a locally finite, connected, infinite graph. 
	Then 
	\be\label{eq: relation of p_T for edge cutset and vertex cutset}
	p_{\mathrm{T,E}}\leq  p_{\mathrm{T,V}}
	\ee
	If moreover $G$ has bounded degree, then the equality holds in \eqref{eq: relation of p_T for edge cutset and vertex cutset}.
\end{lemma}
\begin{proof}[Proof of Lemma \ref{lem: p_T,E less than P_T,V}]
	For \eqref{eq: relation of p_T for edge cutset and vertex cutset}, if $p>p_{\mathrm{T,V}}$, then $\mathbb{E}_p[|C(x)|_V]=\infty$. Since $C(x)$ is connected, $|C(x)|_E\geq |C(x)|_V-1$. Hence $\mathbb{E}_p[|C(x)|_E]=\infty$. Therefore if $p>p_{\mathrm{T,V}}$, then $p\geq p_{\mathrm{T,E}}$. Thus $p_{\mathrm{T,E}}\leq p_{\mathrm{T,V}}$ as desired. 
	
	If $G$ has bounded degree, i.e., $D(G):=\sup\{\mathrm{deg}(v)\colon v\in V \}<\infty$, then by 
	$|C(x)|_E\leq D(G)|C(x)|_V$ one can get the other direction similarly. Hence if $G$ has bounded degree, then $p_{\mathrm{T,V}}(G)=p_{\mathrm{T,E}}(G)$.	
\end{proof}

\begin{example}\label{example: p_T,E less than p_T,V}
	Here we give an example $G$ with unbounded degree and such that $p_{\mathrm{T,E}}< p_{\mathrm{T,V}}$. Let $M>1$ be an integer. Let $C_n$ be a complete graph with $M^n$ vertices. Let $o=(0,0)$ be the origin of $\mathbb{Z}^2$ and let $(n,0)\in\mathbb{Z}^2,n\geq 1$ be the points on the $x$-axis. For each $n\geq 1$, add an edge from $(n,0)$ to each vertex of $C_n$. Let $G$ be the graph obtained in this way; see Figure \ref{fig: p_T,E less than p_T,V}. Then obviously $p_{\mathrm{c}}(G)=p_{\mathrm{c}}(\mathbb{Z}^2)=\frac{1}{2}$. Note that for $p\in(0,p_{\mathrm{c}})$, $\mathbb{P}_p[o\longleftrightarrow (n,0)\textnormal{ in }G]=\mathbb{P}_p[o\longleftrightarrow(n,0)\textnormal{ in }\mathbb{Z}^2]\approx e^{-n\varphi(p)}$, where $\varphi(p)$ is the reciprocal of the correlation length (see Proposition 6.47 in \cite{Grimmett1999} for example.) When computing $\mathbb{E}_p[|C_o|_V]$, each clique $C_n$ contributes roughly $p\cdot e^{-n\varphi(p)}\cdot M^n$ but when computing $\mathbb{E}_p[|C_o|_E]$, each clique $C_n$ contributes roughly $p^2\cdot e^{-n\varphi(p)}\cdot M^{2n}$. Using the properties of $\varphi(p)$ (Theorem 6.14 in \cite{Grimmett1999}) it is easy to show that $0<p_{\mathrm{T,E}}(G)=\varphi^{-1}(2\log M)<p_{\mathrm{T,V}}(G)=\varphi^{-1}(\log M)<p_{\mathrm{c}}(G)$ and we omit the details.
	\begin{figure}[h!]
		\centering
		\begin{tikzpicture}[scale=0.7, text height=1.5ex,text depth=.25ex] 
			%\draw [help lines] (0,0) grid (15,8);
			
			%%%%%%%%%%%%%   draw the lattice Z^2  %%%%%%%%%%%%%%
			
			%%%%%  horizontal lines %%%%%	
			\draw [color=teal] (-1.5,1.5) to (10.5,1.5);
			
			\draw [color=teal] (0,3) to (12,3);
			
			\draw [color=teal] (1.5,4.5) to (13.5,4.5);
			
			\draw [color=teal] (3,6) to (15,6);
			
			%	\draw [color=teal] (-1,7.5) to (15,7.5);

			%%%%%% vertical lines   %%%%%%%%%
			\draw [color=teal] (-1,1) to (5,7);
			
			\draw [color=teal] (1,1) to (7,7);
			
			\draw [color=teal] (3,1) to (9,7);
			
			\draw [color=teal] (5,1) to (11,7);
			
			\draw [color=teal] (7,1) to (13,7);
			
			\draw [color=teal] (9,1) to (15,7);

			%%%%%%%%%%%%    add the cliques   %%%%%%%%%%%%%%%%%

			\draw[fill=black] (3,3) circle [radius=0.06];
			\node[below] at (3,3) {$(0,0)$};

			\draw[fill=black] (5,3) circle [radius=0.06];
			\node[below] at (5,3) {$(1,0)$};
			
			\begin{scope}[shift={(2,0)}]
				\draw [color=red] (3,5) ellipse (0.6 and 0.3);
				\draw [color=red] (3,3) to (2+0.4,5);
				\draw [color=red] (3,3) to (4-0.4,5);
				\node [below,color=red] at  (3,5.2) {$\cdots$};
				\node[above, color=red] at (3,5.2) {$C_1$};
			\end{scope}
			
			\draw[fill=black] (7,3) circle [radius=0.06];
			\node[below] at (7,3) {$(2,0)$};
			
			\begin{scope}[shift={(4,0)}]
				\draw [color=red] (3,5) ellipse (0.6 and 0.3);
				\draw [color=red] (3,3) to (2+0.4,5);
				\draw [color=red] (3,3) to (4-0.4,5);
				\node [below,color=red] at  (3,5.2) {$\cdots$};
				\node[above, color=red] at (3,5.2) {$C_2$};
			\end{scope}

			\draw[fill=black] (9,3) circle [radius=0.06];
			\node[below] at (9,3) {$(3,0)$};
			
			\begin{scope}[shift={(6,0)}]
				\draw [color=red] (3,5) ellipse (0.6 and 0.3);
				\draw [color=red] (3,3) to (2+0.4,5);
				\draw [color=red] (3,3) to (4-0.4,5);
				\node [below,color=red] at  (3,5.2) {$\cdots$};
				\node[above, color=red] at (3,5.2) {$C_3$};
			\end{scope}	
			
			\draw[fill=black] (11,3) circle [radius=0.06];
			\node[below] at (11,3) {$(4,0)$};
			
			\begin{scope}[shift={(8,0)}]
				\draw [color=red] (3,5) ellipse (0.6 and 0.3);
				\draw [color=red] (3,3) to (2+0.4,5);
				\draw [color=red] (3,3) to (4-0.4,5);
				\node [below,color=red] at  (3,5.2) {$\cdots$};
				\node[above, color=red] at (3,5.2) {$C_4$};
			\end{scope}

		\end{tikzpicture}
		\caption{An example with $0<p_{\mathrm{T,E}}<p_{\mathrm{T,V}}<p_{\mathrm{c}}<1$.}
		\label{fig: p_T,E less than p_T,V}
	\end{figure}
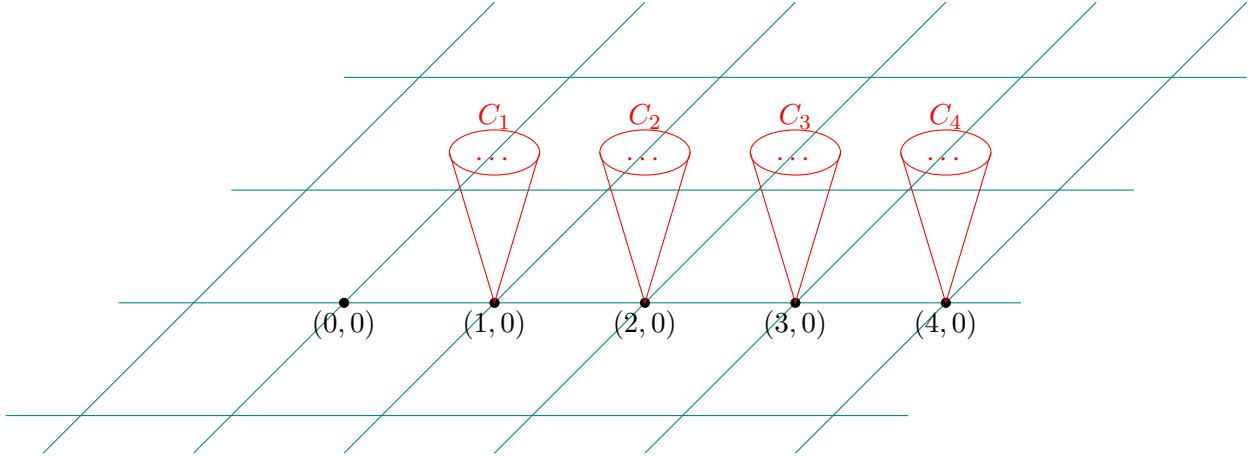

\end{example}

Before proving Lemma \ref{lem: p_cut,e less than p_cut,v}, we recall the definitions of boundaries of a set of vertices. 
\begin{definition}\label{def: vertex boundary and edge boundary}
	For a nonempty set of vertices $K\subset V$, we define its inner vertex boundary, outer vertex boundary and edge boundary as follows. 
	The \notion{inner vertex boundary}	$\partial^{\mathrm{in}}_{\mathrm{V}}K$ is 
	\[
	\partial^{\mathrm{in}}_{\mathrm{V}}K:=\{y\in K\colon \exists\, z\notin K \,\,\mathrm{ s.t. }\,\,y\sim z \},
	\]
	where $y\sim z$ denotes that $y$ and $z$ are neighbors in $G$.
	The \notion{outer vertex boundary}	$\partial_{\mathrm{V}}K$ is 
	\[
	\partial_{\mathrm{V}}K:=\{z\notin K\colon \exists\, y\in K \,\,\mathrm{ s.t. }\,\,y\sim z \}.
	\]
	The \notion{edge boundary} $\Delta K$ is defined as 
	\[
	\Delta K:=\{e=(y,z)\in E\colon\, y\in K,z\notin K \}.
	\]
\end{definition}

\begin{lemma}\label{lem: minimal cutset as boundary}
	If $\Pi$ is a minimal vertex cutset separating $x$ from infinity, then $\partial_V  S(\Pi)=\Pi$, where $S(\Pi)$ is the connected component of $x$ in the subgraph  $G\backslash \Pi$. Similarly, if  $\Pi$ is a minimal edge cutset separating $x$ from infinity, then $\Delta  S(\Pi)=\Pi$, where $S(\Pi)$ is the connected component of $x$ in the subgraph  $G\backslash \Pi$.
\end{lemma}

\begin{definition}\label{def: boundary cutset}
	Suppose $G$ is a locally finite, connected, infinite graph. Fix $x\in V(G)$. Let $\mathscr{B}_E=\mathscr{B}_E(x)$ be the collection of edge cutsets that are  the edge boundary of some finite cluster of $x$.  In other words,
	\[
	\mathscr{B}_E=\{\Delta S\colon  S \textnormal{ is a finite connected subgraph containing }x \}.
	\]
	We call $\mathscr{B}_E$ the family of \notion{boundary edge cutset}.
	Similarly we denote by $\mathscr{B}_V$ the collection of vertex cutsets that arise as outer vertex boundary of some finite cluster of $x$ and call $\mathscr{B}_V$ the family of \notion{boundary vertex cutset}.
\end{definition}

\begin{proof}[Proof of Lemma \ref{lem: minimal cutset as boundary}]
	Suppose $\Pi$ is a minimal vertex cutset separating $x$ from infinity. The connected component $S(\Pi)$ of $x$ in $G\backslash \Pi$ is finite. 
	
	First we show that $\partial_V S(\Pi)\subset \Pi$. For any $z\in \partial_V S(\Pi)$, by definition of $\partial_V S(\Pi)$, there is some vertex  $y\in S(\Pi)$ such that $y\sim z$. Then if $z\notin \Pi$, then by definition of $S(\Pi)$, then $z$ can be connected to $x$ via a path from $x$ to $y$ in $G\backslash \Pi$ and the edge $(y,z)$. This implies that $z\in S(\Pi)$ if $z\notin \Pi$, which contradicts with the choice of $z\in\partial_V S(\Pi)$.  Hence $\partial_V S(\Pi)\subset \Pi$. 
	
	On the other hand, since $\partial_V S(\Pi)$ is also a vertex cutset and $\Pi$ is minimal with respect to inclusion, one has that $\partial_V S(\Pi)\supset \Pi$.
	
	The case of minimal edge cutset can be proved similarly and we omit the details. 
\end{proof}

\begin{remark}
	The reverse of Lemma \ref{lem: minimal cutset as boundary} is not true. For example consider the half integer line $G=(\mathbb{N},E)$, where $E=\{(n,n+1)\colon n\in\mathbb{N}\}$.  Let $S=\{10,11,\ldots,100\}$, then $\Delta S=\{(9,10),(100,101) \}$ is not minimal. 
\end{remark}

\begin{proof}[Proof of Lemma \ref{lem: p_cut,e less than p_cut,v}]

	For  any edge cutset $\Pi_E$ separating $x$ from infinity, let $S(\Pi_E)$ be the connected component of $x$ in $G\backslash \Pi_E$.  Since $\Pi_E$ is a cutset, $S(\Pi_E)$ is finite.  Let $\Pi_V=\Pi_V(\Pi_E)$ be the endpoints of edges in $\Delta S(\Pi_E)$ that are not in $S(\Pi_E)$. Then $\Pi_V$ is a vertex cutset, since every infinite path from $x$ to infinity has to leave $S(\Pi_E)$, the first vertex on the path that is not in $S(\Pi_E)$ must be a vertex in $\Pi_V$.   For each $v\in\Pi_V$, pick an arbitrary edge $e=e(v)\in\Delta S(\Pi_E)$ such that $e$ is incident to $v$. Note that 
	\begin{enumerate}
		\item $\Delta S(\Pi_E)\subset \Pi_E$;
		\item for distinct $v\in\Pi_V$, the edges $e(v)$ are also distinct (since each such edge $e(v)$ has exactly one endpoints not in $S(\Pi_E)$, i.e., $v$);
		\item for such $v$ and $e=e(v)$, $\mathbb{P}[e\in C(x)]\geq p\mathbb{P}_p[v\in C(x)]$ (By insertion-tolerance of Bernoulli bond percolation, see \cite[Exercise 7.1]{LP2016}).
	\end{enumerate}
	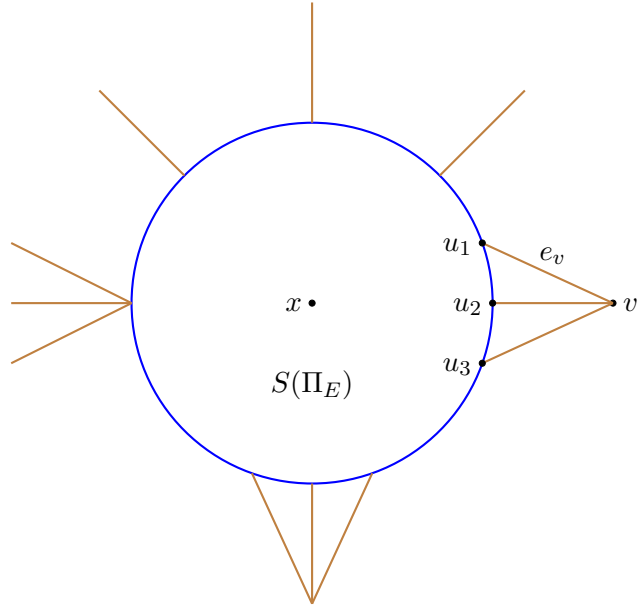
\begin{figure}[h!]
		\centering
		\begin{tikzpicture}[scale=0.7, text height=1.5ex,text depth=.25ex] 
			%\draw [help lines] (0,0) grid (10,10);
			\draw[fill=black] (5,5) circle [radius=0.05];
			\node[left] at (5,5) {$x$};
			\node[below] at (5,4) {$S(\Pi_E)$}; 
			\draw[color=blue, thick] (5,5) circle (3);
			
			\draw[fill=black] (10,5) circle [radius=0.05];
			\node[right] at (10,5) {$v$};
			
			\draw[color=brown, thick] (8,5)--(10,5);
			\draw[color=brown, thick] (7.828,4)--(10,5);
			\draw[color=brown, thick] (7.828,6)--(10,5);
			
			\draw[fill=black] (7.828,4) circle [radius=0.05];
			\node[left] at (7.828,4) {$u_3$};
			\draw[fill=black] (8,5) circle [radius=0.05];
			\node[left] at (8,5) {$u_2$};
			\draw[fill=black] (7.828,6) circle [radius=0.05];
			\node[left] at (7.828,6) {$u_1$};
			
			\node[above] at (9,5.5) {$e_v$};
			%\node[below] at (9,4.5) {$e'$};
			
			\draw[color=brown, thick] (5,2)--(5,0);
			\draw[color=brown, thick] (6,2.172)--(5,0);
			\draw[color=brown, thick] (4,2.172)--(5,0);
			
			\draw[color=brown, thick] (5,8)--(5,10);
			\draw[color=brown, thick] (7.121,7.121)--(8.536,8.536);
			\draw[color=brown, thick] (10-7.121,7.121)--(10-8.536,8.536);

			\draw[color=brown, thick] (2,5)--(0,5);
			\draw[color=brown, thick] (2,5)--(0,6);
			\draw[color=brown, thick] (2,5)--(0,4);
		\end{tikzpicture}
		\caption{A systematic drawing of $S(\Pi_E)$, $v\in\Pi_V(\Pi_E)$ and $e=e(v)$; edges in $\Delta S(\Pi_E)$ are colored brown.}
		\label{fig: edge cutset vs vertex cutset}
	\end{figure}

	Hence 
	\be\label{eq: vertex cutset intersection less than 1/p times of edge cutset}
	\mathbb{E}_p[|C(x)\cap \Pi_V|]=\sum_{v\in\Pi_V}{\mathbb{P}_p[v\in C(x)]}
	\leq \frac{1}{p}\sum_{e\in \Pi_E}{\mathbb{P}_p[e\in C(x)]}=\frac{1}{p}\mathbb{E}_p[|C(x)\cap \Pi_E|]
	\ee
	If $0<p<p_{\mathrm{cut,E}}$, then $\inf_{\Pi_E}\mathbb{E}_p[|C(x)\cap \Pi_E|]=0$. Thus by \eqref{eq: vertex cutset intersection less than 1/p times of edge cutset}, for $p<p_{\mathrm{cut,E}}$
	\[
	\inf_{\Pi_V}\mathbb{E}_p[|C(x)\cap \Pi_V|]=0.
	\]
	Hence if $p<p_{\mathrm{cut,E}}$, then $p\leq p_{\mathrm{cut,V}}$. This implies \eqref{eq: relation of p_cut for edge cutset and vertex cutset}.

	Now we assume that $G$ has bounded degree. Let $\Pi_V$ be a vertex cutset and without loss of generality we assume that $x\notin \Pi_V$. Let $S(\Pi_V)$ denote the  connected component of $x$ in $G\backslash \Pi_V$. Since $\Pi_V$ is a cutset, the connected component of $S(\Pi_V)$ is finite. Let $\Pi_E=\Pi_E(\Pi_V)$ be the edge boundary $\Delta S(\Pi_V)$. Now for each edge $e\in\Delta S(\Pi_V)$, there is a unique vertex  $v=v(e)\in\Pi_V$ associated to $e$:  $v$ is incident to $e$ (the other endpoint of $e$ is in $S(\Pi_V)$, which is disjoint from $\Pi_V$ by its definition). Note that 
	\begin{enumerate}
		\item for each $v\in\Pi_V$, there are at most $D=D(G)$ edges in $\Delta S(\Pi_V)$ associated to it; 
		\item for each $e\in \Delta S(\Pi_V)$ and its associated vertex $v=v(E)\in \Pi_V$, $\mathbb{P}_p[v\in C(x)]\geq \mathbb{P}_p[e\in C(x)]$. 
	\end{enumerate} 
	Hence 
	\be\label{eq: edge cutset intersection less than D times of vertex cutset}
	\mathbb{E}_p[|C(x)\cap \Pi_E|]=\sum_{e\in\Pi_E}{\mathbb{P}_p[e\in C(x)]}
	\leq D\sum_{v\in \Pi_V}{\mathbb{P}_p[v\in C(x)]}=D\mathbb{E}_p[|C(x)\cap \Pi_V|]
	\ee
	Thus when $D<\infty$,  by \eqref{eq: edge cutset intersection less than D times of vertex cutset} one has that
	\[
	\forall\,\,p<p_{\mathrm{cut,V}},\,\,\,\inf_{\Pi_E}\mathbb{E}_p[|C(x)\cap \Pi_E|]=0.
	\]
	Hence when $D<\infty$, $p<p_{\mathrm{cut,V}}\Rightarrow p\leq p_{\mathrm{cut,E}}$. Together with \eqref{eq: relation of p_cut for edge cutset and vertex cutset} one has  the equality when $D<\infty$.	
\end{proof}

When considering $p_{\mathrm{cut,V}}$ and $p_{\mathrm{cut,E}}$, obviously it suffices to consider minimal cutsets (with respect to inclusion). However a priori  it is not clear whether it is sufficient to consider minimal cutsets for $p'_{\mathrm{cut,E}}$ and $p'_{\mathrm{cut,V}}$.  This consideration for minimal cutsets together with Lemma \ref{lem: minimal cutset as boundary} motivates us to consider $p''_{\mathrm{cut,E}}$ and $p''_{\mathrm{cut,V}}$ in Definition \ref{def: def of p''_cut} and it turns out they all coincides with each other; see Theorem \ref{thm: general version}.
\begin{definition}\label{def: def of p''_cut}
	Suppose $G$ is a locally finite, connected, infinite graph. Fix $x\in V(G)$.
	Define 
	\begin{equation*}\label{eq: def of p''_cut,V in bond percolation}
		p''_{\mathrm{cut,V}}= p''_{\mathrm{cut,V}}(G):=\sup\Big\{p\geq0\colon \inf_{\Pi\in\mathscr{B}_V}\sum_{v\in \Pi}\mathbb{P}_p[A(x,v,\Pi)]=0 \Big\},
	\end{equation*}
	where the infimum is taken over all \textbf{boundary vertex cutsets}  $\Pi$ that  separate $x$ from infinity.
	
	Similarly, we define 
	\begin{equation*}\label{eq: def of p''_cut,E in bond percolation}
		p''_{\mathrm{cut,E}}= p''_{\mathrm{cut,E}}(G):=\sup\Big\{p\geq0\colon \inf_{\Pi_E\in\mathscr{B}_E}\sum_{e\in \Pi_E}\mathbb{P}_p[A(x,e,\Pi_E)]=0 \Big\},
	\end{equation*}
	where the infimum is taken over all \textbf{boundary edge cutsets} $\Pi_E$ that separate $x$ from infinity. 
	
\end{definition}

By Definition \ref{def: def of p'_cut,E and p'_cut V} and \ref{def: def of p''_cut},  and inequalities \eqref{eq: p'_cut,E leq p_c}, \eqref{eq: p'_cut,V leq p_c} one has that 
\be\label{eq: p''_cut,E leq p'_cut,E leq p_c}
p''_{\mathrm{cut,E}}\leq p'_{\mathrm{cut,E}}\leq p_{\mathrm{c}}\,\,\,\mathrm{and}\,\,\,p''_{\mathrm{cut,V}}\leq p'_{\mathrm{cut,V}}\leq p_c.
\ee
Theorem \ref{thm: p_c=p_cut'} is contained in  the following more general theorem.
\begin{theorem}\label{thm: general version}
	For Bernoulli bond percolation on every locally finite, connected, infinite graph $G$, one has that 
	\[
	p''_{\mathrm{cut,E}}=p'_{\mathrm{cut,E}}=p''_{\mathrm{cut,V}}=	p'_{\mathrm{cut,V}}=p_{\mathrm{c}}.
	\]
\end{theorem}
\begin{lemma}\label{lem: p''_cut,E less than p''_cut,V}
	Suppose $G$ is a locally finite, connected infinite graph.  Then 
	\be\label{eq: p'_cut,E less than p'_cut,V}
	p''_{\mathrm{cut,E}}\leq p''_{\mathrm{cut,V}}.
	\ee
\end{lemma}
\begin{proof}
	Let $\Pi_E \in \mathscr{B}_E$ be a boundary edge cutset separating $x$ from infinity. Let $S$ be the finite connected component of $x$ in $G\backslash \Pi_E$. By definition $\Pi_E=\Delta S$.  Let $\Pi_V=\partial_V S$ be the outer vertex boundary of $S$.  Then $\Pi_V \in \mathscr{B}_V$ is a boundary vertex cutset separating $x$ from infinity. 
	
	For each $v\in\Pi_V$,  if the event $A(x,v,\Pi_V)$ occurs, then there is a self-avoiding open path $\gamma_{x,v}$ from $x$ to $v$ only using $v$ in $\Pi_V$. Hence this path uses only one edge $e$ in $\Delta S$, namely the edge $e$ on $\gamma_{x,v}$ that is incident to $v$. Hence the event $A(x,e,\Pi_E)$ occurs for this  edge $e$ on the path $\gamma_{x,v}$.  Thus 
	\[
	\mathbb{P}_p(A(x,v,\Pi_V))\leq \sum_{e\in\Delta S\colon e\sim v} \mathbb{P}_p(A(x,e,\Pi_E)),
	\]
	where $e\sim v$ denotes that $v$ is an endpoint of $e$.
	
	Note that for any two distinct vertices $v,v'\in \Pi_V$, the two sets $\{e\in\Delta S\colon e\sim v\}$ and $\{e'\in\Delta S\colon e'\sim v' \}$ are disjoint. Hence summing the above inequality over $v\in \Pi_V=\partial S$, 
	one has that 
	\[
	\sum_{v\in \Pi_V}\mathbb{P}_p(A(x,v,\Pi_V))\leq \sum_{e\in \Pi_E}\mathbb{P}_p(A(x,e,\Pi_E)).
	\]
	From this we have that  $p<p''_{\mathrm{cut,E}}\Rightarrow p\leq p''_{\mathrm{cut,V}}$ and then we have the desired inequality \eqref{eq: p'_cut,E less than p'_cut,V}. 
\end{proof}

\section{Proof of Theorem \ref{thm: general version}}\label{sec: proof of the main result}

Duminil-Copin and Tassion \cite{DCT15} gave a new proof of the sharpness of the phase transition \cite{AB1987,Men1986}. For our purpose, we just need to look at the short version for Bernoulli percolation \cite{DCT15short}. 

For $p\in[0,1]$, $x\in V$ and a finite set $S$ with $x\in S\subset V$,  define 
\be\label{def: definition of phi p}
\varphi_p(x,S):=p\sum_{y\in S}\sum_{z\notin S,(y,z)\in E}\mathbb{P}_p(x\stackrel{S}{\longleftrightarrow} y),
\ee
where $\{x\stackrel{S}{\longleftrightarrow} y \}$ denotes the event that there is an open path connecting $x$ and $y$ only using vertices lying in $S$. Recall that the edge boundary $\Delta S$ of $S$ is the set of edges that connect $S$ to its complement. So $\varphi_p(x,S)$ is the expected number of open edges on the edge boundary $\Delta S$ which has an endpoint is connected to $x$ via an open path entirely lying in $S$. For transitive graphs, Duminil-Copin and Tassion defined
\[
\widetilde{p}_c:=\sup\{p\geq0\colon \inf\{\varphi_p(x,S)\colon x\in S, \, S \textnormal{ is finite}\} <1 \}
\] 
and showed that $\widetilde{p}_c=p_{\mathrm{c}}$ for transitive graphs.

The main new ingredient of the proof of Theorem \ref{thm: general version} is the following observation.
\begin{proposition}\label{prop: observation}
	For Bernoulli bond percolation on a locally finite, connected, infinite graph $G$ one has that
	\be\label{eq: key observation}
	\inf_{\Pi\in\mathscr{B}_E} \sum_{e\in \Pi}\mathbb{P}_p(A(x,e,\Pi))=\inf_S\varphi_p(x,S),
	\ee
	where the infimum on the left hand side of \eqref{eq: key observation} is over all the \textbf{boundary edge cutsets} separating $x$ from infinity and the infimum on the right is over all finite sets containing $x$. 
\end{proposition}
We have that  $p_{\mathrm{c}}(G)=\sup\{p\geq0\colon \inf\{\varphi_p(x,S)\colon x\in S, \, S \textnormal{ is finite}\} =0 \}$ for all locally finite, connected, infinite graphs in light of Proposition \ref{prop: observation} and Theorem \ref{thm: general version}.

\begin{proof}[Proof of Proposition \ref{prop: observation}]
	On the one hand, for any finite set $S$ containing $x$, let $S'$ be the connected component of $x$ in the induced subgraph of $S$. Then $\Pi(S):=\Delta S'$ is a boundary edge cutset separating $x$ from infinity. For each edge $e=(y,z)\in \Delta S'$, say $y\in S',z\notin S'$, it is easy to see that $z\notin S$ and 
	\[
	\mathbb{P}_p(A(x,e,\Pi(S)))= p\cdot \mathbb{P}_p[x\stackrel{S}{\longleftrightarrow} y].
	\]
	%% Actually both sides are zero if $e\in \Delta S \backslash \Delta S(x)$, where $S(x)$ deonotes those vertices in S with a path to x that lying entirely in S. %%
	Summing this over all edges $e\in\Delta S'$, one has that
	\[
	\inf_{\Pi\in\mathscr{B}_E} \sum_{e\in \Pi}\mathbb{P}_p(A(x,e,\Pi))\leq \sum_{e\in \Pi(S)}\mathbb{P}_p(A(x,e,\Pi(S)))= \varphi_p(x,S')=\varphi_p(x,S),
	\]
	where the last equality is a simple observation from the definition of $S'$.
	
	Hence 
	\[
	\inf_{\Pi\in\mathscr{B}_E} \sum_{e\in \Pi}\mathbb{P}_p(A(x,e,\Pi))\leq \inf_S\varphi_p(x,S).
	\]

	On the other hand, for any boundary edge cutset $\Pi$ separating $x$ from infinity, let $S=S(\Pi)$ be the connected component of $x$ in the graph $G\backslash \Pi$. By Definition \ref{def: boundary cutset},  $\Delta S= \Pi$.
	For each edge $e=(y,z)\in \Delta S= \Pi$, say $y\in S,z\notin S$, one has that 
	\[
	\mathbb{P}_p(A(x,e,\Pi))=p\cdot \mathbb{P}_p[x\stackrel{S}{\longleftrightarrow} y].
	\]
	Summing this over all edges $e\in\Delta S$, one has that for a boundary edge cutset $\Pi$ and $S=S(\Pi)$
	\[
	\sum_{e\in \Pi}\mathbb{P}_p(A(x,e,\Pi))= \varphi_p(x,S(\Pi))\geq \inf_S\varphi_p(x,S).
	\]
	Hence one has the other direction
	\[
	\inf_{\Pi\in\mathscr{B}_E}  \sum_{e\in \Pi}\mathbb{P}_p(A(x,e,\Pi))\geq \inf_S\varphi_p(x,S).\qedhere
	\]
\end{proof}

Next we recall a lemma from \cite{DCT15short}. For a finite set $\Lambda$, let $\Lambda^c$ denote its complement in $V$. Let $\Lambda_n$ denote the ball $\{y: d(x,y)\leq n\}$ of radius $n$ centered at $x$, where $d$ denotes the graph distance on $G$. 
\begin{lemma}[Lemma 2.1 of \cite{DCT15short}]\label{lem: differential inequality}
	For $x\in V$ and ball $\Lambda_n$ with $n\geq1$, one has 
	\be\label{eq: differential ineq}
	\frac{d }{dp}\mathbb{P}_p(x\longleftrightarrow \Lambda_n^c)\geq \frac{1}{p(1-p)}\cdot \inf_{S\subset \Lambda_n,x\in S}\varphi_p(x,S)\cdot [1-\mathbb{P}_p(x\longleftrightarrow \Lambda_n^c)]
	\ee
	
\end{lemma}

\begin{proof}[Proof of Theorem \ref{thm: general version}]
	By \eqref{eq: p''_cut,E leq p'_cut,E leq p_c} and Lemma \ref{lem: p''_cut,E less than p''_cut,V}, it suffices to show  $p''_{\mathrm{cut,E}}\geq p_{\mathrm{c}}$. 
	
	Suppose $p''_{\mathrm{cut,E}}<p_{\mathrm{c}}$. Pick $p_0,p_1$ such that $	p''_{\mathrm{cut,E}}<p_0<p_1<p_{\mathrm{c}}$.
	
	By the definition of $	p''_{\mathrm{cut,E}}$ and Proposition \ref{prop: observation}, there is a constant $\kappa>0$ such that for any $p\in[p_0,p_1]$, 
	\[
	\inf_S\varphi_p(x,S)\geq \kappa.
	\]
	Write $\theta_x(p):=\mathbb{P}_p(x\longleftrightarrow\infty)$ and $\theta_{x,n}(p):=\mathbb{P}_p(x\longleftrightarrow\Lambda_n^c)$. 
	By \eqref{eq: differential ineq} one has that for $p\in[p_0,p_1]$, 
	\[
	\frac{\theta_{x,n}'(p)}{1-\theta_{x,n}(p)}\geq \frac{\kappa}{p(1-p)}.
	\]
	Integrating this inequality from $p_0$ to $p_1$, one has that 
	\be\label{eq: integration of the differential inequality}
	\theta_{x,n}(p_1)\geq 1-\left(\frac{1-p_1}{p_1}\cdot\frac{p_0}{1-p_0}\right)^{\kappa}+\theta_{x,n}(p_0)\left(\frac{1-p_1}{p_1}\cdot\frac{p_0}{1-p_0}\right)^{\kappa}\geq 1-\left(\frac{1-p_1}{p_1}\cdot\frac{p_0}{1-p_0}\right)^{\kappa}.
	\ee
	Letting $n \rightarrow \infty$ one has that 
	\[
	\theta_x(p_1)\geq 1-\left(\frac{1-p_1}{p_1}\cdot\frac{p_0}{1-p_0}\right)^{\kappa}>0,
	\]
	which contradicts with the choice that $p_1<p_{\mathrm{c}}$. Hence $p''_{\mathrm{cut,E}}\geq p_{\mathrm{c}}$ and we are done.
\end{proof}

\section{Percolation probability for subperiodic trees}\label{sec: periodic trees}
To highlight the importance of Proposition \ref{prop: observation}, in this section we discuss some applications of it to subperiodic trees.

%\subsection{Remarks and questions on subperiodic trees}
For transitive graphs, Duminil-Copin and Tassion pointed out that $\inf\{\varphi_p(x,S): x\in S, \, S \textnormal{ is finite}\} \geq1$  at $p=p_{\mathrm{c}}$. Using this they obtained a lower bound for percolation probability on transitive graphs with $p_{\mathrm{c}}\in(0,1)$: $\theta(p)\geq \frac{p-p_{\mathrm{c}}}{p(1-p_{\mathrm{c}})}$ for $p\geq p_{\mathrm{c}}$. Here for transitive graphs, the percolation probability $\theta(x,p)$ does not depend on $x$ and we simply write it as $\theta(p)$. This lower bound can be extended to 0-subperiodic trees. We first adopt some notations and then recall the definitions of periodic and subperiodic trees as in \cite[Section 3.3]{LP2016}. 

\noindent\textbf{Notation.} Suppose $T$ is an infinite, locally finite tree with a distinguished vertex $o$, called the \notion{root} of $T$. Write $|x|$ for the graph distance from $x$ to $o$; $x\leq y$ if $x$ is on the shortest path from $o$ to $y$; $x<y$ if $x\leq y$ and $x\neq y$; $x\rightarrow y$ if $x\leq y$ and $|y|=|x|+1$ and in this case we call $x$ the parent of $y$; and $T^x$ for the subtree of $T$ containing the vertices $y\geq x$. For Bernoulli$(p)$ percolation on the tree $T$ with root $o$, let $\theta(p):=\mathbb{P}_p[o\longleftrightarrow\infty]$ be the probability that $o$ is in an infinite cluster. 

\begin{definition}[Definition on page 82 of \cite{LP2016}]
	Let $N\geq0$ be an integer. An infinite, locally finite tree $T$ with root $x$ is called \notion{$N$-periodic} (resp., \notion{$N$-subperiodic}), if $\forall x\in T$ there exists an adjacency-preserving bijection (resp., injection) $f:T^x\rightarrow T^{f(x)}$ with $|f(x)|\leq N$. A tree is \notion{periodic} (resp., \notion{subperiodic}) if there is some $N\geq0$ for which it is $N$-periodic (resp., $N$-subperiodic).
\end{definition}

\begin{remark}
	For a 0-subperiodic tree $T$ with root $o$, one has that $\inf_{\Pi}\sum_{e\in\Pi}\textnormal{br}(T)^{-|e|}\geq 1$ (formula (3.7) on page 85 of \cite{LP2016}). It is well-known that  $p_{\mathrm{c}}(T)=1/\textnormal{br}(T)$ for every locally finite infinite tree $T$ (for instance see \cite[Theorem 5.15]{LP2016}). Hence 
	\[
	\inf_{\Pi\in\mathscr{B}_E} \sum_{e\in\Pi} \mathbb{P}_{p_{\mathrm{c}}}(A(o,e,\Pi))=\inf_{\Pi\in\mathscr{B}_E}\sum_{e\in\Pi}p_{\mathrm{c}}^{|e|}\geq\inf_{\Pi}\sum_{e\in\Pi}\textnormal{br}(T)^{-|e|}\geq 1.
	\]
	Then by Proposition \ref{prop: observation} one can set $p_0=p_{\mathrm{c}}$ and $\kappa=1$ in \eqref{eq: integration of the differential inequality} and letting $n\rightarrow\infty$ to get 
	\[
	\theta(p)=\mathbb{P}_p(o\longleftrightarrow\infty)\geq \frac{p-p_{\mathrm{c}}}{p(1-p_{\mathrm{c}})},\,\, p\geq p_{\mathrm{c}}
	\]
	for every 0-subperiodic tree $T$.
\end{remark}
Actually this is true for all subperiodic trees with nontrivial $p_c$. 
\begin{proposition}\label{prop: positive lower right Dini derivative}
	Consider Bernoulli percolation on a subperiodic tree $T$ with $p_{\mathrm{c}}(T)<1$.
	Then the lower right Dini derivative of the percolation probability $\theta(p)$ at $p_{\mathrm{c}}$ belong to $(0,\infty]$.
\end{proposition}	
\begin{proof}
	Theorem 3.8 in \cite{LP2016} says that $\inf_{\Pi}\sum_{e\in\Pi}\textnormal{br}(T)^{-|e|}>0$ for a general subperiodic tree $T$ with $p_{\mathrm{c}}(T)<1$, where $\Pi$ runs over all edge cutsets separating the root of $T$ from infinity. 
	Define \[
	\alpha(o,p):=\inf\{\varphi_p(o,S): o\in S, \, S \textnormal{ is finite}\}.
	\]
	Then as before, one has that for a subperiodic tree $T$ with root $o$ and $p_{\mathrm{c}}<1$, \[
	\alpha(o,p_{\mathrm{c}})\stackrel{\textnormal{Prop.}\ref{prop: observation}}{=}	\inf_{\Pi\in\mathscr{B}_E} \sum_{e\in\Pi} \mathbb{P}_{p_{\mathrm{c}}}(A(o,e,\Pi))=\inf_{\Pi\in\mathscr{B}_E}\sum_{e\in\Pi}p_{\mathrm{c}}^{|e|}\geq\inf_{\Pi}\sum_{e\in\Pi}\textnormal{br}(T)^{-|e|}>0.
	\]
	By the definition of $\varphi$ in \eqref{def: definition of phi p} it is obvious that $\alpha(o,p)$ is increasing in $p$. 
	Setting $p_0=p_{\mathrm{c}}$ in \eqref{eq: integration of the differential inequality} and letting $n\rightarrow\infty$ one has that 
	\[
	\frac{\theta(p)-\theta(p_{\mathrm{c}})}{1-\theta(p_{\mathrm{c}})}\geq 1- \left(\frac{1-p}{p}\cdot \frac{p_{\mathrm{c}}}{1-p_{\mathrm{c}}}\right)^{\alpha(o,p_{\mathrm{c}})}. 
	\]
	This implies that the lower right Dini derivative of the percolation probability $\theta(p)$ at $p_{\mathrm{c}}$ is positive:
	\[
	D_+\theta(p_{\mathrm{c}}):=\liminf_{p\rightarrow p_{\mathrm{c}}^+} \frac{\theta(p)-\theta(p_{\mathrm{c}})}{p-p_{\mathrm{c}}}\geq 
	\frac{\alpha(o,p_{\mathrm{c}})(1-\theta(p_{\mathrm{c}}))}{p_{\mathrm{c}}(1-p_{\mathrm{c}})}
	>0.
	\] 
\end{proof}
Without of the restriction of (sub)periodicity, it is possible to have the lower right Dini derivative of the percolation probability being equal to zero. See the following two remarks. 
\begin{remark}\label{rem: a spherically symmetric tree}
	It is easy to construct trees with the property that at $p=p_{\mathrm{c}}$,
	\[
	\inf\{\varphi_p(o,S): o\in S, \, S \textnormal{ is finite}\}=0.
	\]
	Indeed, we construct a spherically symmetric tree $T$ with root $o$ as follows. Let $T_n$ denote the set of vertices with graph distance $n$ to the root $o$. If $n=2^k$ for some $k\geq0$, let each vertex in $T_n$ have exactly one child; otherwise, let each vertex in $T_n$ have exactly two children. Then it is easy to see that 
	\[
	|T_n|\asymp \frac{2^n}{n}. 
	\]
	Here for two positive function $f,g$ on $\mathbb{Z}^+$, $f(n)\asymp g(n)$ means that there exist constants $c_1,c_2>0$ such that $c_1g(n)\leq f(n)\leq c_2g(n)$ for all $n>0$. 
	
	Hence $\textnormal{br}(T)=\liminf_n|T_n|^{1/n}=2$ by Exercise 1.2 in \cite{LP2016}. Thus $p_{\mathrm{c}}(T)=1/\textnormal{br}(T)=1/2$. Let $S_n$ be the ball of radius $n$ and center $x$. Then at $p=p_{\mathrm{c}}=1/2$,
	\[
	\varphi_p(o,S_n)=|T_{n+1}|\cdot \frac{1}{2^{n+1}}\leq \frac{c_2}{n+1}.
	\]
	Thus $\inf\{\varphi_p(o,S): o\in S, \, S \textnormal{ is finite}\}=0.$
\end{remark}

\begin{remark}
	For the spherically symmetric tree $T$ in Remark \ref{rem: a spherically symmetric tree}, one also has 
	\[
	\theta(p)\asymp \big(p-\frac{1}{2}\big)^2,\,\,p\geq p_{\mathrm{c}}
	\]
	and in particular, the upper right Dini derivative $D^+\theta(p_{\mathrm{c}}):=\limsup_{p\rightarrow p_{\mathrm{c}}^+} \frac{\theta(p)-\theta(p_{\mathrm{c}})}{p-p_{\mathrm{c}}}=0$. 
	
	In fact, let $c(e)=(1-p)^{-1}p^{|e|}$ be the conductance of edge $e$. Then formula (5.12) on page 142 of \cite{LP2016} is satisfied with $\mathbf{P}=\mathbb{P}_p$. 
	Since $T$ is spherically symmetric, for $p>p_{\mathrm{c}}=1/2$,   the effective resistance is
	\[
	\mathscr{R}(o\longleftrightarrow \infty)=\sum_{n=1}^{\infty}(1-p)p^{-n}/|T_n|\asymp\sum_{n=1}^{\infty}
	\frac{(1-p)n}{(2p)^n}\asymp \frac{1-p}{(2p-1)^2}. 
	\]
	Then by Theorem 5.24 \cite{LP2016} one has 
	\[
	\theta(p)\asymp\frac{\mathscr{C}(o\longleftrightarrow \infty)}{1+\mathscr{C}(o\longleftrightarrow \infty)}=\frac{1}{1+\mathscr{R}(o\longleftrightarrow \infty)}\asymp\big(p-\frac{1}{2}\big)^2,\,\,p>p_{\mathrm{c}}.
	\]
\end{remark}

Proposition \ref{prop: positive lower right Dini derivative} states that the lower right Dini derivative of the percolation probability on a subperiodic tree with  $p_c\in(0,1)$ is positive at $p_c$ and it might equal to infinity in some cases (Example \ref{example: binary tree and Fibonacci tree and an infinite right derivative example}). This leads us to the following question:
\begin{question} \label{ques: q2}
	What kind of subperiodic trees have the property that the right Dini derivatives of $\theta(p)$ at $p_{\mathrm{c}}$ are finite? What kind of subperiodic trees have the property that $\lim_{p\downarrow p_{\mathrm{c}}} \frac{\theta(p)-\theta(p_{\mathrm{c}})}{p-p_{\mathrm{c}}}\in(0,\infty)$?
\end{question}

The critical exponent $\beta$ for Bernoulli percolation is characterized by $\theta(p)-\theta(p_{\mathrm{c}})\approx (p-p_{\mathrm{c}})^\beta$. For $\mathbb{Z}^2$, it is conjectured that $\theta(p)-\theta(p_{\mathrm{c}})\approx (p-p_{\mathrm{c}})^\beta$ for  $\beta=\frac{5}{36}$ \cite[Table 10.1 on page 279]{Grimmett1999}, in particular in this case the lower right Dini derivative at $p_{\mathrm{c}}$ is infinite \cite{KestenZhang1987}. Indeed for site percolation on the triangular lattice in the plane, one does have $\theta(p)-\theta(p_{\mathrm{c}})= (p-p_{\mathrm{c}})^{\frac{5}{36}+o(1)}$ \cite[Theorem 1.1]{SmirnovWerner2001}.  Question \ref{ques: q2} asks  what kind of subperiodic trees have $\beta=1$.

A partial answer for  Question \ref{ques: q2} is Theorem \ref{thm: finite upper right Dini derivative for certain periodic tree} which considers directed covers of strongly connected graphs.  An oriented graph $G$ is called \notion{strongly connected} if for any two vertices $u,v$ of $G$, there is a directed path in $G$ from $u$ to $v$. 
Suppose that $G$ is a finite oriented graph and $v$ is any vertex in $G$. The \notion{directed cover} of $G$ based at $v$ is the tree $T$ whose vertices are the finite paths of edges $\langle e_1,e_2,\ldots,e_n\rangle$ in
$G$ that start at $v$. We take the root of $T$ to be the empty path and we join two vertices in $T$ by an edge when a path is an extension of the other path by one more edge in $G$. Every periodic tree is a directed cover of a finite directed graph $G$; for a proof see pages 82-83 of \cite{LP2016}. Also not all periodic trees have finite right derivatives for $\theta(p)$ at $p_{\mathrm{c}}$; see item 3 in Example \ref{example: binary tree and Fibonacci tree and an infinite right derivative example}. 
%\begin{lemma}\label{lem: periodic trees as directed covers}
%	Every periodic tree is a directed cover of a finite directed graph $G$.  
%\end{lemma}
\begin{theorem}\label{thm: finite upper right Dini derivative for certain periodic tree}
	If $T$ is a periodic tree with root $o$ and $p_{\mathrm{c}}(T)\in(0,1)$ and  it is the directed cover of some strongly connected graph, then the right derivative of $\theta(p)$ exists at $p_{\mathrm{c}}$ and this derivative is positive and  finite. 
\end{theorem}

We first outline the ideas of the proof of Theorem \ref{thm: finite upper right Dini derivative for certain periodic tree}. Suppose $T$ is the directed cover of some finite strongly connected directed graph $G=(V,E)$ based at some vertex $v_1\in V$, where  $V=\{v_1,\ldots,v_n\}$. Let $\theta_i(p)$ be the percolation probability of the tree $T_i$ that is the directed cover of $G$ based at $v_i$. Then these quantities $\theta_i(p)$ are related by a family of algebraic equations \eqref{eq: relations of percolation prob via adjacency matrix}.  We will then proceed as follows:
\begin{enumerate}
	\item[1.] We first show that  $\theta_i(p)$ are continuous at $p_c(T_i)$ (Lemma \ref{lem: theta(p_c) equals zero for periodic trees}).
	
	\item[2.] Then we show that these trees $T_i$ have the same critical probability $p_c$  and the upper right Dini derivatives of $\theta_i(p)$ are finite at $p_c$ (Lemma \ref{lem: finite upper right Dini derivative}).
	
	\item[3.] We finish the preparation by showing  that the right derivatives of $\theta_i(p)$ at $p_c$ exist via a convexity/concavity argument in Lemma \ref{lem: remain concave or convex near p_c}. The proof of Lemma \ref{lem: remain concave or convex near p_c} is a little bit involved and we will need two more lemmas for its proof:
	\begin{itemize}
		\item the functions $\theta_i(p)$ are analytic on $(p_c,1)$ (Lemma \ref{lem: infinite differentiability}) and 
		\item the uniqueness of $\theta_i(p)$ as solutions of \eqref{eq: relations of percolation prob via adjacency matrix} (Lemma \ref{lem: uniquess as solution of relations of perc prob}).
	\end{itemize} 
	
	\item[4.] We then finish the proof of Theorem \ref{thm: finite upper right Dini derivative for certain periodic tree}: the existence is from Step 3 (Lemma \ref{lem: remain concave or convex near p_c}), the positiveness is from Proposition \ref{prop: positive lower right Dini derivative} and the finiteness is from Step 2 (Lemma \ref{lem: finite upper right Dini derivative}). 
\end{enumerate}

\begin{lemma}\label{lem: theta(p_c) equals zero for periodic trees}
	Suppose $T$ is a periodic tree with root $o$ and $p_{\mathrm{c}}(T)<1$. Then 
	$\theta(p_{\mathrm{c}})=0$.
\end{lemma}
\begin{proof}
	At $p=p_{\mathrm{c}}=\frac{1}{\mathrm{br}(T)}$, if we put conductance $c(e(x))=(1-p_{\mathrm{c}})^{-1}p_{\mathrm{c}}^{|x|}$, where $e(x)$ is the edge from $x$ to its parent, then (5.12) on page 142 of \cite{LP2016} is
	satisfied. As noted on  page 142 line 17 of  \cite{LP2016}, these conductances correspond to the homesick random walk $\mathrm{RW}_{\mathrm{br}(T)}$. If we put resistance $\Phi(e(x))=\lambda^{|x|-1}$ for the edge $e(x)$ instead, then it is known that as $\lambda\uparrow \lambda_*=\frac{1}{p_c}$, the effective resistance  from the root to infinity of the corresponding network is tend to infinity  \cite[Theorem 5.1]{Lyons1990}. This implies that  the homesick random walk $\mathrm{RW}_{\mathrm{br}(T)}$ is recurrent. Hence by Corollary 5.25 of \cite{LP2016} we know $\theta(p_{\mathrm{c}})=0$.
\end{proof}

Now we restrict to a subset of periodic trees that are directed covers of finite strongly connected oriented graphs.
\begin{lemma}\label{lem: finite upper right Dini derivative}
	Suppose $G=(V(G),E(G))$ is a finite, strongly connected directed graph and $V(G)=\{v_1,\ldots,v_n\}$. Let $\lambda_*$ be the largest eigenvalue of the adjacency matrix $A_G$ of $G$.  Let $T_i$ be the directed cover of $G$ based at $v_i$ and denote its root by $o_i$.  Then $p_{\mathrm{c}}(T_1)=\cdots=p_{\mathrm{c}}(T_n)=\frac{1}{\lambda_*}$. 
	
	Moreover if $\lambda_*>1$, then the upper right Dini derivative of $\theta_i(p)$ at $p_{\mathrm{c}}$ is finite for every $i\in\{ 1,\cdots,n \}$, where $\theta_i(p):=\mathbb{P}_p[o_i\longleftrightarrow \infty \textnormal{ in }T_i]$ denotes the probability that the root $o_i$ of $T_i$ is in an infinite open cluster. 
\end{lemma}
\begin{proof}
	The first part is a standard result. See the discussion on pages 83-84 of \cite{LP2016} for example. 
	
	Since $G$ is strongly connected, $A_G$ is irreducible. Hence by the Perron-Frobenius theorem (e.g.\ see \cite[Theorem 8.4.4]{HJ2013matrixanalysis}) there is  a left $\lambda_*$-eigenvector $v_*=(v_1,\ldots,v_n)$ all of whose entries are positive.  We also normalize $v_*$ such that its $l_2$-norm  is $1$.

	Since  $o_i\not\leftrightarrow\infty$ in  $T_i$ if and only if $o_i$ can't connect to infinity via any of its children, we have the following relations for these percolation probabilities:
	\[
	1-\theta_i(p)=\prod_{j=1}^{n}{[1-p\theta_j(p)]^{a_{ij}}},\,i\in\{1,\ldots,n\},
	\]
	i.e.,
	\be\label{eq: relations of percolation prob via adjacency matrix}
	\theta_i(p)=1-\prod_{j=1}^{n}{[1-p\theta_j(p)]^{a_{ij}}},\,i\in\{1,\ldots,n\},
	\ee
	where $a_{ij}$ is the $(i,j)$-entry of the matrix $A_G$, i.e., the number of directed edges in $G$ from vertex $v_i$ to $v_j$.

	Denote by $\theta_{\mathrm{max}}(p)=\max\{\theta_1(p),\cdots,\theta_n(p)\}$. Since $G$ is strongly connected, there exists $M>0$ such that there is a directed path with length at most $M$ from $v_i$ to $v_j$ for any pair $v_i,v_j\in V(G)$. Hence $\theta_i(p)\geq p^M\theta_{\mathrm{max}}(p)$ for all $i=1,\ldots,n$. Thus for $p> p_{\mathrm{c}}$,  
	\be\label{eq: same order for perc. prob.}
	0<\theta_i(p)\asymp \theta_{\mathrm{max}}(p),\,i\in\{1,\ldots,n\},. 
	\ee
	
	By Lemma \ref{lem: theta(p_c) equals zero for periodic trees} and the right continuity of $\theta_i(p)$ (e.g., see \cite[Exercise 7.33]{LP2016}), one has that 
	\be\label{eq: theta_i(p) is little o(p-p_c)}
	0<\theta_i(p)=o(p-p_{\mathrm{c}}),\,\,0< p-p_{\mathrm{c}}\ll1.
	\ee
	Using \eqref{eq: theta_i(p) is little o(p-p_c)} when $0<p-p_{\mathrm{c}}\ll1$ and $i\in\{1,\ldots,n\},$ we can rewrite \eqref{eq: relations of percolation prob via adjacency matrix} as 
	\be\label{eq: rewritten relations of perc. prob.}
	\theta_i(p)=p\sum_{j=1}^{n}a_{ij}\theta_j(p)-p^2\sum_{j=1}^{n}\binom{a_{ij}}{2}\theta_j^2(p)-p^2\sum_{j\neq k}{a_{ij}a_{ik}\theta_j(p)\theta_k(p)}+\theta_{\mathrm{max}}^2(p)\cdot o(1),
	\ee
	where we use the convention that $\binom{a_{ij}}{2}=0$ if $a_{ij}=0,1$.
	
	Multiplying $v_i$ on both sides of \eqref{eq: rewritten relations of perc. prob.} and adding them up, one has that 
	\begin{eqnarray}\label{eq: averaged perc prob via left eigenvector}
		\sum_{i=1}^{n}v_i\theta_i(p)&=&p\sum_{i=1}^{n}v_i\sum_{j=1}^{n}a_{ij}\theta_j(p)+\theta_{\mathrm{max}}^2(p)\cdot o(1)\nonumber\\
		&&-p^2	\sum_{i=1}^{n}v_i\bigg[\sum_{j=1}^{n}\binom{a_{ij}}{2}\theta_j^2(p)+\sum_{j\neq k}{a_{ij}a_{ik}\theta_j(p)\theta_k(p)}\bigg]
	\end{eqnarray}

	Since $p_{\mathrm{c}}(T_i)=\frac{1}{\lambda_*}<1$, there exists some $i$ such that either $a_{ij}\geq 2$ for some $j$ or $a_{ij}a_{ik}\geq 1$ for some $j\neq k$.  Therefore by \eqref{eq: same order for perc. prob.} and \eqref{eq: averaged perc prob via left eigenvector} there exists $c>0$ such that 
	\be\label{eq: inequality for theta_max's square}
	\sum_{i=1}^{n}v_i\theta_i(p)\leq p\sum_{i=1}^{n}v_i\sum_{j=1}^{n}a_{ij}\theta_j(p)-cp^2\theta_{\mathrm{max}}^2(p),\,\,0< p-p_{\mathrm{c}}\ll1.
	\ee
	
	Since $v_*$ is a left $\lambda_*$-eigenvector of $A_G$, one has that 
	\be\label{eq: eigenvector equation}
	\sum_{i=1}^{n}v_i\sum_{j=1}^{n}a_{ij}\theta_j(p)=v_*A_G\bm{\theta}(p)=\lambda_*v_*\cdot \bm{\theta}(p)=\lambda_*\sum_{i=1}^{n}v_i\theta_i(p),
	\ee
	where $\bm{\theta}(p)=(\theta_1(p),\cdots,\theta_n(p))^T$ is the vector of percolation probabilities in $\mathbb{R}^n$.
	
	Plugging \eqref{eq: eigenvector equation} into \eqref{eq: inequality for theta_max's square} and using $\lambda_*=\frac{1}{p_{\mathrm{c}}}$  and \eqref{eq: same order for perc. prob.} one has that for $0<p-p_{\mathrm{c}}\ll 1$, 
	\[
	cp^2\theta_{\mathrm{max}}^2(p)\leq \frac{p-p_{\mathrm{c}}}{p_{\mathrm{c}}}\sum_{i=1}^{n}v_i\theta_i(p)\leq c'\theta_{\mathrm{max}}(p)(p-p_{\mathrm{c}}),
	\]
	for some constant $c'>0$. This implies that $\theta_i(p)\leq \theta_{\mathrm{max}}(p)\leq c''(p-p_{\mathrm{c}})$ for $0<p-p_{\mathrm{c}}\ll1$ and then we have the desired result on the upper right Dini derivative. 
\end{proof}
The strongly connectedness of $G$ is needed for the finiteness of the right Dini derivative. See the following example. 

\begin{figure}[h!]
	\centering
	\begin{tikzpicture}[scale=.8, text height=1.5ex,text depth=.25ex] 
		%\draw [help lines] (0,0) grid (15,5);

		%%%%%%%%%%%%%  G_a is just two loops at a vertex %%%%%%%%%%%%%
		\draw[fill=black] (2,2) circle [radius=0.06];
		\node[left] at (2,2) {$v_1$};
		
		\node [below] at (2,0) {$G_a$};
		
		\begin{scope}[shift={(1,2)}]
			\tikzset{style={->},>=stealth,-{Latex[length=2mm,width=1mm]}}
			\tkzDefPoints{0/0/O}
			\tkzDrawArc[R, color=blue](O,1)(0,359.8)
		\end{scope}

		\begin{scope}[shift={(3,2)}]
			\tikzset{style={->},>=stealth,-{Latex[length=2mm,width=1mm]}}
			\tkzDefPoints{0/0/O}
			\tkzDrawArc[R, color=blue](O,1)(180.1,179.9)
		\end{scope}

		%%%%%%%%%%%%%%% G_b %%%%%%%%%%%%%%
		
		\draw[fill=black] (7,2) circle [radius=0.06];
		\node[left] at (7,2) {$v_1$};
		
		\draw[fill=black] (9,2) circle [radius=0.06];
		\node[right] at (9,2) {$v_2$};
		
		\node[below] at (7,0)  {$G_b$};
		
		\begin{scope}[shift={(6,2)}]
			\tikzset{style={->},>=stealth,-{Latex[length=2mm,width=1mm]}}
			\tkzDefPoints{0/0/O}
			\tkzDrawArc[R, color=blue](O,1)(0,359.8)
		\end{scope}
		
		\draw[->,-{Latex[length=2mm,width=1mm]},>=stealth,color=blue] (7,2) to[bend right] (9,2);
		
		\draw[->,-{Latex[length=2mm,width=1mm]},>=stealth,color=blue] (9,2) to[bend right] (7,2);

		%%%%%%%%%%%%%%%%%  G_c %%%%%%%%%%%%%%%%%%
		
		\begin{scope}[shift={(6,0.5)}]
			
			\draw[fill=black] (7,3) circle [radius=0.06];
			\node[left] at (7,3) {$v_1$};
			
			\draw[fill=black] (9,3) circle [radius=0.06];
			\node[right] at (9,3) {$v_2$};

			\begin{scope}[shift={(6,3)}]
				\tikzset{style={->},>=stealth,-{Latex[length=2mm,width=1mm]}}
				\tkzDefPoints{0/0/O}
				\tkzDrawArc[R, color=blue](O,1)(0,359.8)
			\end{scope}
			
			\draw[->,-{Latex[length=2mm,width=1mm]},>=stealth,color=blue] (7,3) to[bend right] (9,3);
			
			\draw[->,-{Latex[length=2mm,width=1mm]},>=stealth,color=blue] (9,3) to[bend right] (7,3);
			
		\end{scope}

		\begin{scope}[shift={(6,-2)}]
			
			\draw[fill=black] (7,3) circle [radius=0.06];
			\node[left] at (7,3) {$v_4$};
			
			\draw[fill=black] (9,3) circle [radius=0.06];
			\node[right] at (9,3) {$v_3$};

			\begin{scope}[shift={(6,3)}]
				\tikzset{style={->},>=stealth,-{Latex[length=2mm,width=1mm]}}
				\tkzDefPoints{0/0/O}
				\tkzDrawArc[R, color=blue](O,1)(0,359.8)
			\end{scope}
			
			\draw[->,-{Latex[length=2mm,width=1mm]},>=stealth,color=blue] (7,3) to[bend right] (9,3);
			
			\draw[->,-{Latex[length=2mm,width=1mm]},>=stealth,color=blue] (9,3) to[bend right] (7,3);
			
		\end{scope}
		
		\draw[->,-{Latex[length=2.5mm,width=1mm]},>=stealth,color=blue] (15,3.5) to (15,1);
		
		\node[below] at (13,0) {$G_c$};
		
	\end{tikzpicture}
	\caption{Directed graphs $G_a,G_b,G_c$ from left to right.}
	\label{fig: three examples}
\end{figure}
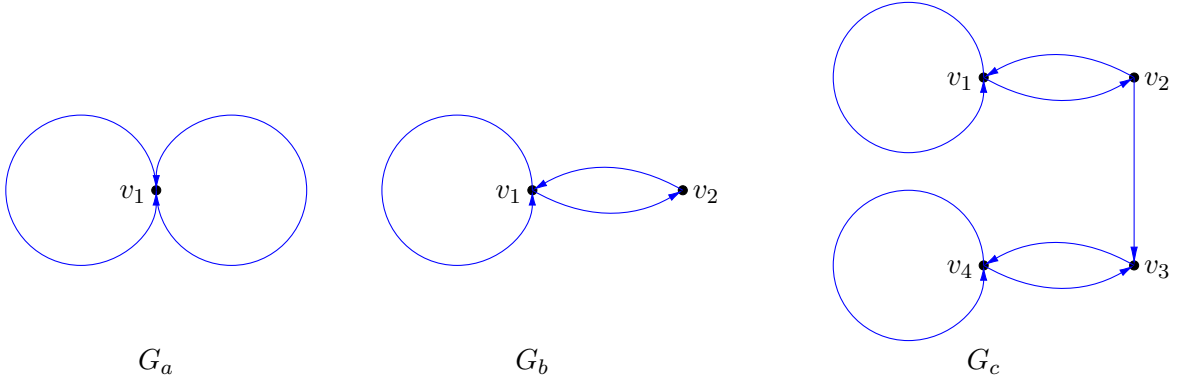

\begin{example}\label{example: binary tree and Fibonacci tree and an infinite right derivative example}
	Let $G_a,G_b,G_c$ be as illustrated in Figure \ref{fig: three examples}. Let $T_s,s\in\{a,b,c\}$ be the directed cover of $G_s$ based at $v_1(G_s)$.  
	\begin{enumerate}
		\item The tree $T_a$ is a binary tree with root $o$. It is easy to see that $p_{\mathrm{c}}(T_a)=\frac{1}{2}$ and for $p\geq \frac{1}{2}$, $\theta(p)=\frac{2p-1}{p^2}$. In this case $\theta'_+(p_{\mathrm{c}})=8$ and $\theta(p)$ is concave on $(p_{\mathrm{c}},1)$.
		
		\item The tree $T_b$ is a Fibonacci tree with root $o$ and $\mathrm{deg}(o)=2$. See Figure 3.2 on page 83 of \cite{LP2016} for an illustration of the Fibonacci tree. It is easy to see that $p_{\mathrm{c}}(T_b)=\frac{\sqrt{5}-1}{2}$.  Writing $\theta(p)=\mathbb{P}_p[o\longleftrightarrow\infty]$ and using
		\eqref{eq: relations of percolation prob via adjacency matrix} one has that $\theta(p)=\frac{p^2+p-1}{p^3}$  for $p\geq \frac{\sqrt{5}-1}{2}$.  Hence in this case $\theta'_+(p_{\mathrm{c}})=5+\sqrt{5}$ and $\theta(p)$ is also concave on $(p_{\mathrm{c}},1)$.
		
		\item The tree $T_c$ also has $p_{\mathrm{c}}(T_c)=\frac{\sqrt{5}-1}{2}$. Actually if we define $T_i,\theta_i(p)$ as in Lemma \ref{lem: finite upper right Dini derivative}, then it is easy to check that $p_{\mathrm{c}}(T_i)=\frac{\sqrt{5}-1}{2},i=1,2,3,4$. Solving \eqref{eq: relations of percolation prob via adjacency matrix} one can find that for $p\in(p_{\mathrm{c}},1)$,
		\be
		\left\{
		\begin{array}{ccc}
			\theta_1(p)&=&\frac{(1-2p)(p^2+p-1)+\sqrt{(p^2+p-1)(-3p^2+5p-1)}}{2p^2(1-p^2)}\\
			&&\\
			\theta_2(p)&=&\frac{p^2+p-1+\sqrt{(p^2+p-1)(-3p^2+5p-1)}}{2p^2}\\
			&&\\
			\theta_3(p)&=&\frac{p^2+p-1}{p^2}\\
			&&\\
			\theta_4(p)&=&\frac{p^2+p-1}{p^3}
		\end{array}
		\right.
		\ee
		In particular, $\theta(p)=\theta_1(p)=\Theta(\sqrt{p-p_{\mathrm{c}}})$ for $0<p-p_{\mathrm{c}}\ll1$ and thus the right Dini derivative at $p_{\mathrm{c}}$ is infinite. 
		One can also check that $\theta_1(p)$ and $\theta_2(p)$ are concave on $(p_{\mathrm{c}},1)$. 
	\end{enumerate}
\end{example}
Example \ref{example: binary tree and Fibonacci tree and an infinite right derivative example} and the fact that $\theta(p)-\theta(p_{\mathrm{c}})\approx (p-p_{\mathrm{c}})^{5/36}$ on the triangle lattice on the plane \cite[Theorem 1.1]{SmirnovWerner2001} lead us to the following question:
\begin{question}
	For a transitive graph or a periodic tree with root $o$, is the percolation probability $\theta(p)$ concave on $(p_{\mathrm{c}},1)$?
\end{question}

We now proceed to show the right derivatives of $\theta_i(p)$ exist at $p_c$.
\begin{lemma}\label{lem: remain concave or convex near p_c}
	The right derivative of $\theta_i(p)$ exists at $p_c$ for all $i\in\{1,\ldots,n\}$. 
\end{lemma}

As mentioned earlier, the proof of Lemma \ref{lem: remain concave or convex near p_c} is somewhat long and we begin by showing that $\theta_i(p)$ is analytic on $(p_c,1)$ for all $i\in\{1,\ldots,n\}$. 

\begin{lemma}\label{lem: infinite differentiability}
	Under the same assumptions as Lemma \ref{lem: finite upper right Dini derivative} one has that  $\theta_i(p)$ is analytic on $(p_{\mathrm{c}},1)$ for all $i\in\{1,\cdots,n \}$. 
\end{lemma}

\begin{proof}
	Recall that the percolation probabilities satisfy  \eqref{eq: relations of percolation prob via adjacency matrix}:
	\[
	\theta_i(p)=1-\prod_{j=1}^{n}{[1-p\theta_j(p)]^{a_{ij}}},\,i\in\{1,\ldots,n\},
	\]
	where $a_{ij}$ is  the number of directed edges in $G$ from vertex $v_i$ to $v_j$.
	
	Define $f_i:\mathbb{R}^{1+n}\rightarrow\mathbb{R}$ for $i\in\{1,\ldots,n\}$ by 
	\[
	f_i((x,y_1,\ldots,y_n)^T)=y_i-1+\prod_{j=1}^{n}[1-xy_j]^{a_{ij}}. 
	\]
	Write $\bm{f}=(f_1,\cdots,f_n)^T$. 
	By  \eqref{eq: relations of percolation prob via adjacency matrix}, we know that $(p,\theta_1(p),\ldots,\theta_n(p))$ is a positive solution of $\bm{f}=\bm{0}$ when $p>p_{\mathrm{c}}$. 
	
	Note that when $i\neq j$,
	\[
	\frac{\partial f_i}{\partial y_j}(p,\theta_1(p),\ldots,\theta_n(p))=
	-pa_{ij}[1-p\theta_j(p)]^{a_{ij}-1}\cdot \prod_{j'\neq j}{[1-p\theta_{j'}(p)]^{a_{ij'}}}\stackrel{\eqref{eq: relations of percolation prob via adjacency matrix}}{=}\frac{-pa_{ij}(1-\theta_i(p))}{1-p\theta_j(p)}
	\]
	and 
	\[
	\frac{\partial f_i}{\partial y_i}(p,\theta_1(p),\ldots,\theta_n(p))=1	-pa_{ii}[1-p\theta_i(p)]^{a_{ii}-1}\cdot \prod_{j'\neq i}{[1-p\theta_{j'}(p)]^{a_{ij'}}}\stackrel{\eqref{eq: relations of percolation prob via adjacency matrix}}{=}1-\frac{pa_{ii}(1-\theta_i(p))}{1-p\theta_i(p)}.
	\]
	Therefore the Jacobi matrix  $J=\big[\frac{\partial f_i}{\partial y_j}(p,\theta_1(p),\ldots,\theta_n(p))\big]_{1\leq i,j\leq n}$ can be written as 
	\be\label{eq: Jacobi matrix}
	J=I-BC
	\ee
	where $I$ is the identity matrix and $B$ is a diagonal matrix with $b_{ii}=1-\theta_i(p)$ and $C$ is a matrix with $(i,j)$-entry $c_{ij}=\frac{pa_{ij}}{1-p\theta_j(p)}$.

	%$
	%A=\begin{bmatrix}
	%	1-\theta_1(p)   &               &             &    \\
	%	                & 1-\theta_2(p) &             &     \\
	%	                &               &  \cdots           &      \\
	%	                &               &             & 1-\theta_n(p)
	%\end{bmatrix},
	%B=\begin{bmatrix}
	%   &               &             &    \\
	%&                      &             &     \\
	%&               &        \frac{pa_{ij}}{1-p\theta_j(p)}    &      \\
	%&               &             & 
	%\end{bmatrix}$

	Notice that Bernoulli$(p)$ percolation on $T_i$ can also be viewed as a multi-type Galton--Watson tree $Z$. Each vertex $u$ on the tree corresponds to a directed path on $G$. If the endpoint of the path is $v_j$, then we say that $u$ has type $j$. In particular, we view the root of $T_i$ is of type $i$. The number of type $j$ children of a type $i$ vertex has Binomial distribution $\mathrm{Bin}(a_{ij},p)$. The percolation probability $\theta_i(p)$ is just the non-extinction probability for such a $n$-type  Galton--Watson tree started with a single type $i$ vertex. Let $\mathbb{P}_s$ and $\mathbb{E}_s$ denote the probability measure and corresponding expectation for such an  $n$-type Galton-Watson tree started with a single ancestor with type $s\in\{ 1,\ldots,n \}$.
	
	Now  let $\mathsf{Ext}$ denote the event that   the $n$-type Galton--Watson tree is extinct. Then $\mathbb{P}_i[\mathsf{Ext}]=1-\theta_i(p)$. Let $Z_{1j}$ denote the number of children of type $j$ of $Z_0$. For a nonnegative integer sequence $(t_1,\ldots,t_n)$ with $t_i\leq a_{ij}$, one has that 
	\begin{eqnarray}\label{eq: distribution conditioned on extinction}
		&&\mathbb{P}_i\big[  Z_{1j}=t_j,j=1,2,\ldots,n   \big|  \mathsf{Ext}\big]\nonumber\\
		&=&\frac{1}{1-\theta_i(p)}\cdot \prod_{j=1}^{n}\binom{a_{ij}}{t_j}p^{t_j}(1-p)^{a_{ij}-t_j}\cdot (1-\theta_j(p))^{t_j}\nonumber\\
		&=&\frac{1}{1-\theta_i(p)}\cdot \prod_{j=1}^{n}\binom{a_{ij}}{t_j}(1-p)^{a_{ij}-t_j}\cdot (p-p\theta_j(p))^{t_j}
		\nonumber\\
		&\stackrel{\eqref{eq: relations of percolation prob via adjacency matrix}}{=}&
		\prod_{j=1}^{n}\binom{a_{ij}}{t_j}\Big(\frac{1-p}{1-p\theta_j(p)}\Big)^{a_{ij}-t_j}\cdot \Big(\frac{p-p\theta_j(p)}{1-p\theta_j(p)}\Big)^{t_j}
		\nonumber\\
		&=&\prod_{j=1}^{n}\binom{a_{ij}}{t_j} \Big(\frac{p-p\theta_j(p)}{1-p\theta_j(p)}\Big)^{t_j}\cdot \Big(1-\frac{p-p\theta_j(p)}{1-p\theta_j(p)}\Big)^{a_{ij}-t_j}.
	\end{eqnarray}
	
	By \cite{JL2008branching} we know conditioned on extinction, the $n$-type Galton-Watson tree is still a multi-type Galton--Watson tree.  Let $\widetilde{\mathbb{P}}_s$ and $\widetilde{\mathbb{E}}_s$ denote the probability measure and corresponding expectation for the $n$-type Galton-Watson tree started with a single ancestor with type $s$ \textbf{conditioned on extinction}. By \eqref{eq: distribution conditioned on extinction}, conditioned on extinction, the number of type $j$ children of a type $i$ vertex has Binomial distribution $\mathrm{Bin}(a_{ij},\frac{p-p\theta_j(p)}{1-p\theta_j(p)})$. Hence the mean offspring  matrix $M$ has $(i,j)$-entry $m_{ij}=a_{ij}\frac{p-p\theta_j(p)}{1-p\theta_j(p)}=(1-\theta_j(p))\cdot \frac{pa_{ij}}{1-p\theta_j(p)}$. Observe that 
	\be\label{eq: observation for the mean size matrix conditioned on extinction}
	M=CB
	\ee
	Let $q:=\max_{1\leq j\leq n}[1-\theta_j(p)]$ be the maximum of the extinction probability. For $p>p_{\mathrm{c}}$, we know $q<1$.  Let $Z_k$ denote the size of $k$-th generation of the multi-type Galton-Watson tree. As the last displayed inequality on page 547 of \cite{JL2008branching}, one has that 
	\[
	\widetilde{\mathbb{E}}_s[Z_k]\leq \frac{1}{1-\theta_s(p)}\cdot \mathbb{E}_s\big[Z_kq^{Z_k}\big]\rightarrow 0
	\textnormal{ as }k\rightarrow\infty.
	\]
	Hence the largest eigenvalue $\lambda_1(M)$ for the mean offspring  matrix $M$ satisfies $\lambda_1(M)<1$. By \cite[Theorem 1.3.22]{HJ2013matrixanalysis}, the largest eigenvalue of $BC$ satisfies that $\lambda_1(BC)=\lambda_1(CB)=\lambda_1(M)<1$. Therefore by \eqref{eq: Jacobi matrix} the Jacobi matrix $J$ is invertible for $p\in(p_{\mathrm{c}},1)$. Hence by the analytic implicit function theorem, we obtain that the functions $\theta_i(p)$ are analytic on $(p_{\mathrm{c}},1)$.
\end{proof}

\begin{lemma}\label{lem: uniquess as solution of relations of perc prob}
	The solution $(p,\theta_1(p),\ldots,\theta_n(p))$ of \eqref{eq: relations of percolation prob via adjacency matrix} in $(p_{\mathrm{c}},1)\times(0,1)^n$  is unique.
\end{lemma}
The following Proposition \ref{prop: monotone for the operator B_p} will be needed for Lemma \ref{lem: uniquess as solution of relations of perc prob}.

For  $p\in[0,1]$, we define the operator $\bm{B}_p:[0,1]^n\longrightarrow[0,1]^n$ as given by \eqref{eq: relations of percolation prob via adjacency matrix}:
\be\label{eq: definition of the operator B_p}
\bm{B}_p(\bm{\alpha})_i=1-\prod_{j=1}^{n}[1-p\alpha_j]^{a_{ij}},
\ee
where  $\bm{\alpha}=(\alpha_1,\cdots,\alpha_n)^T\in [0,1]^n$. For example, $\bm{B}_p(\bm{0})=\bm{0}$. 

For $\bm{\alpha},\bm{\beta}\in[0,1]^n$, write $\bm{\alpha}\leq \bm{\beta}$ if $\alpha_i\leq \beta_i$ for all $i\in\{1,\cdots,n\}$ and write $\bm{\alpha}\prec\bm{\beta}$ if $\bm{\alpha}\leq \bm{\beta}$ and $\bm{\alpha}\neq \bm{\beta}$. 

\begin{proposition}\label{prop: monotone for the operator B_p}
	We have  the following properties for the operator $\bm{B}_p$.
	\begin{enumerate}
		\item[(a)] The operator $\bm{B}_p$ is increasing in the sense that if $\bm{\alpha}\leq \bm{\beta}$, then $\bm{B}_p(\bm{\alpha})\leq \bm{B}_p(\bm{\beta})$. 
		
		\item[(b)] Moreover, if $p\in(0,1)$, then $\bm{B}_p$ is strictly increasing in the sense that if $\bm{\alpha}\prec \bm{\beta}$, then $\bm{B}_p(\bm{\alpha})\prec \bm{B}_p(\bm{\beta})$.
		
		\item[(c)] If $0\leq p_1<p_2\leq 1$ and $\bm{\alpha}\in[0,1]^n$, then $\bm{B}_{p_1}(\bm{\alpha})\leq \bm{B}_{p_2}(\bm{\alpha})$. 
		
		\item[(d)] Moreover, if $\bm{\alpha}\in[0,1]^n$ and $\bm{\alpha}\neq \bm{0}$, then for $0\leq p_1<p_2< 1$, one has that $\bm{B}_{p_1}(\bm{\alpha})\prec \bm{B}_{p_2}(\bm{\alpha})$. 
		
		\item[(e)] For $p>0$, if $\bm{0}\neq \bm{\alpha}\in[0,1]^n$ is a fixed point of $\bm{B}_p$, i.e., $\bm{B}_p(\bm{\alpha})=\bm{\alpha}$, then $\alpha_i>0$ for  all $i\in \{ 1,\cdots,n \}$.
	\end{enumerate} 
\end{proposition}
\begin{proof}[Proof of Proposition \ref{prop: monotone for the operator B_p}]
	The items (a) and (c) are obvious from the definition of $\bm{B}_p$. For item (b), suppose $\alpha_j<\beta_j$ for some $j\in\{ 1,\cdots,n \}$. Since $G$ is strongly connected, there exists some $i$ such that $a_{ij}\geq1$. Then 
	\[
	[1-p\alpha_j]^{a_{ij}}>[1-p\beta_j]^{a_{ij}}
	\]
	and $[1-p\alpha_{j'}]^{a_{ij'}}\geq [1-p\beta_{j'}]^{a_{ij'}}\geq[1-p]^{a_{ij'}} >0$ for $j'\neq j$.
	Therefore $\bm{B}_p(\bm{\alpha})_i<\bm{B}_p(\bm{\beta})_i$. Together with item (a) we know $\bm{B}_p(\bm{\alpha})\prec \bm{B}_p(\bm{\beta})$.
	
	For item (d), the proof is similar to item (b) and we omit it.
	
	For item (e), if $a_{ij}\geq1$, by  \eqref{eq: definition of the operator B_p} and the fact that $\bm{\alpha}$ is a fixed point, 
	\[
	\alpha_i=\bm{B}_p(\bm{\alpha})_i\geq 1-[1-p\alpha_j]^{a_{ij}}\geq 1-[1-p\alpha_j]=p\alpha_j.
	\]
	Repeating this argument, we get 
	\[
	\alpha_i=\bm{B}_p(\bm{\alpha})_i\geq p^M\alpha_{j'},\forall\, j'\in \{ 1,\cdots,n \},
	\]
	where $M$ is the maximum of the lengths of the shortest oriented paths connecting two points in $G$. Therefore since  $\alpha_j>0$ for some $j$, then  $\alpha_i>0$ for all $i\in \{ 1,\cdots,n \}$.
\end{proof}

\begin{proof}[Proof of Lemma \ref{lem: uniquess as solution of relations of perc prob}]	
	Write $\theta_{i,k}(p):=\mathbb{P}_p[o_i\textnormal{ is connected to level }k \textnormal{ of }T_i]$, where by level $k$ we mean the set of vertices in $T_i$ with graph distance $k$ to the root. Let $\bm{\theta}_{k}(p)=(\theta_{1,k}(p),\cdots,\theta_{n,k}(p))^T\in[0,1]^n$. In particular, $\bm{\theta}_{0}(p)=\bm{1}$. Then as \eqref{eq: relations of percolation prob via adjacency matrix}, one has that 
	\[
	\bm{\theta}_{k+1}(p)=\bm{B}_p(\bm{\theta}_{k}(p))
	\]
	and thus $\bm{\theta}_{k}(p)=\bm{B}_p^{\circ k}(\bm{1})$. By the definition of $\theta_{i,k}(p),\theta_i(p)$,
	\[
	\bm{\theta}(p)=\lim_{k\rightarrow\infty}\bm{\theta}_{k}(p)=\lim_{k\rightarrow\infty}\bm{B}_p^{\circ k}(\bm{1}).
	\]
	Suppose $\bm{\alpha}\in[0,1]^n$ is some fixed point of $\bm{B}_p$, i.e., $\bm{B_p}(\bm{\alpha})=\bm{\alpha}$. Then $\bm{\alpha}\leq \bm{1}$. By item (a) of Proposition \ref{prop: monotone for the operator B_p}, one has that $\bm{\alpha}=\bm{B_p}(\bm{\alpha})\leq\bm{B}_p(\bm{1})=\bm{\theta}_1(p)$, and then $\bm{\alpha}=\bm{B_p}(\bm{\alpha})\leq\bm{B}_p(\bm{\theta}_1(p))=\bm{\theta}_2(p),\cdots$
	In the end, we have \[
	\bm{\alpha}\leq \lim_{k\rightarrow\infty}\bm{\theta}_{k}(p)=\bm{\theta}(p),
	\] 
	i.e., $\bm{\theta}(p)$ is the largest fixed point for the operator $\bm{B}_p$ in $[0,1]^n$.
	
	Now suppose $p\in(p_{\mathrm{c}},1)$ and we have some solution $\bm{\alpha}\in[0,1]^n\backslash\{\bm{0} \}$ for \eqref{eq: relations of percolation prob via adjacency matrix}, i.e., $\bm{\alpha}$ is a nonzero fixed point of $\bm{B}_p$ in $[0,1]^n$. Since $\bm{\alpha}\in[0,1]^n$ is a fixed point of $\bm{B}_p$, we have showed that $\bm{\alpha}\leq \bm{\theta}(p)$. Since $\bm{\theta}(p)\geq \bm{\alpha}\neq \bm{0}$, by item (e) in Proposition \ref{prop: monotone for the operator B_p}, $\bm{\alpha}\in(0,1)^n$. Define $p_1:=\sup\{t\geq p_{\mathrm{c}}\colon \bm{\theta}(t)\leq \bm{\alpha}\}$. 
	
	Since  $\bm{\alpha}\leq \bm{\theta}(p)$ and $\lim_{t\downarrow p_{\mathrm{c}}}\bm{\theta}(t)\rightarrow\bm{0}$ (Lemma \ref{lem: theta(p_c) equals zero for periodic trees}), one has that $p_1\in(p_{\mathrm{c}},p]$.  
	Since $\bm{\theta}(p)$ is infinite differentiable (Lemma \ref{lem: infinite differentiability}) and increasing in $(p_{\mathrm{c}},1)$, one has that $\bm{\theta}(p_1)\leq \bm{\alpha}$ and for some $i\in\{1,\ldots,n \}$, $\theta_i(p_1)=\alpha_i$.  
	
	Since $\bm{\theta}(p_1)$ is a fixed point of $\bm{B}_{p_1}$ and $\bm{\alpha}$ is a fixed point of $\bm{B}_p$, one has that 
	\be\label{eq: compare theta(p_1) with alpha}
	\theta_i(p_1)=1-\prod_{j=1}^{n}[1-p_1\theta_j(p_1)]^{a_{ij}}=
	\alpha_i=1-\prod_{j=1}^{n}[1-p\alpha_j]^{a_{ij}}.
	\ee
	Since $G$ is strongly connected, there is some $j$ such that $a_{ij}\geq1$.  As in the proof of item (e) in Proposition \ref{prop: monotone for the operator B_p},  to satisfy \eqref{eq: compare theta(p_1) with alpha}, by $p_1\leq p$ and $\bm{\theta}(p_1)\leq \bm{\alpha}$ one must have 
	\[
	1-p_1\theta_j(p_1)=1-p\alpha_j.
	\]
	Since $0<\theta_j(p_1)\leq \alpha_j\leq \theta_j(p)$ and $p_{\mathrm{c}}<p_1\leq p<1$, one must have $p_1=p$ and $\alpha_j=\theta_j(p_1)$. Therefore by $p_1=p$ and the continuity of $\bm{\theta}$ one has the other direction $\bm{\alpha}\geq \lim_{t\uparrow p_1}\bm{\theta}(t)=\bm{\theta}(p_1)=\bm{\theta}(p)$. Hence $\alpha=\bm{\theta}(p)$ is the unique solution of \eqref{eq: relations of percolation prob via adjacency matrix} in $(0,1)^n$ for $p\in(p_{\mathrm{c}},1)$. 	
\end{proof}

Now we are ready to prove Lemma \ref{lem: remain concave or convex near p_c}.
\begin{proof}[Proof of Lemma \ref{lem: remain concave or convex near p_c}]
	First by \eqref{eq: relations of percolation prob via adjacency matrix} and Lemma \ref{lem: uniquess as solution of relations of perc prob} we know that the set 
	\[
	\{(p,\theta_1(p),\ldots,\theta_n(p))\in (p_{\mathrm{c}},1)\times(0,1)^n\}
	\]
	is semi-algebraic (Definition 2.1.4 of \cite{BCR1998algebraic_geometry}). By Theorem 2.2.1 of \cite{BCR1998algebraic_geometry}, its projection  $S:=\{(\theta_1(p),\ldots,\theta_n(p))\colon p\in(p_{\mathrm{c}},1)\}$ is also semi-algebraic set and by Definition 2.2.5 the map $p\mapsto(\theta_1(p),\ldots,\theta_n(p))$ is a semi-algebraic map from $(p_{\mathrm{c}},1)$ to $S$ in view of \eqref{eq: relations of percolation prob via adjacency matrix}.
	
	By Theorem 2.2.1 and Proposition 2.2.6 of \cite{BCR1998algebraic_geometry}, the maps $p\mapsto \theta_i(p)$ are also semi-algebraic functions on $(p_{\mathrm{c}},1)$. 
	
	By Lemma \ref{lem: infinite differentiability} we already know that the  functions $p\mapsto \theta_i(p)$ are infinitely differentiable on $(p_{\mathrm{c}},1)$. Hence by Proposition 2.9.1 of \cite{BCR1998algebraic_geometry}, we know the second derivatives $\theta_i''(p)$ are also semi-algebraic functions on $(p_{\mathrm{c}},1)$, i.e., the sets $\{(p,\theta_i''(p))\colon p\in(p_{\mathrm{c}},1)  \}$ are semi-algebraic sets for all $i\in\{1,\ldots,n\}$.

	Hence for every $i\in\{1,\ldots,n\}$,  $\{(p,\theta_i''(p))\colon p\in(p_{\mathrm{c}},1),\theta_i''(p)=0  \}$ is a semi-algebraic set since it is the intersection  of two semi-algebraic sets, $\{(p,\theta_i''(p))\colon p\in(p_{\mathrm{c}},1)  \}$ and $\{(p,0)\colon p\in(p_{\mathrm{c}},1)  \}$.  Thus by Theorem 2.2.1 of \cite{BCR1998algebraic_geometry} for every $i\in\{1,\ldots,n\}$, the projection $\{ p\colon p\in(p_{\mathrm{c}},1),\theta_i''(p)=0 \}$ is a semi-algebraic set. 
	By Proposition 2.1.7 of \cite{BCR1998algebraic_geometry}, for every $i\in\{1,\ldots,n\}$, the set $\{ p\colon p\in(p_{\mathrm{c}},1),\theta_i''(p)=0 \}$ is a \textbf{finite} union of points and open intervals. Hence there exists some $\varepsilon_i>0$ such that $\theta_i''(p)$ cannot change its sign  on $(p_{\mathrm{c}},p_{\mathrm{c}}+\varepsilon_i)$. Taking $\varepsilon=\min\{\varepsilon_i\colon i=1,2,\ldots,n\}$ we have that the function $\theta_i''(p)$ does not change its sign (i.e.\ remains nonnegative or nonpositive) on $(p_{\mathrm{c}},p_{\mathrm{c}}+\varepsilon)$ for all $i\in\{1,2,\ldots,n\}$. 
	
	By Lemma \ref{lem: remain concave or convex near p_c} and the right continuity of the functions $\theta_i(p)$ at $p_{\mathrm{c}}$, we know these functions $\theta_i(p)$ are either convex or concave on $[p_{\mathrm{c}},p_{\mathrm{c}}+\varepsilon)$. Hence  $\lim_{p\downarrow p_{\mathrm{c}}}\frac{\theta_i(p)-\theta_i(p_{\mathrm{c}})}{p-p_{\mathrm{c}}}=\lim_{p\downarrow p_{\mathrm{c}}}\frac{\theta_i(p)-0}{p-p_{\mathrm{c}}}$ exists, i.e., the right derivative of $\theta_i(p)$ at $p_{\mathrm{c}}$ exists. 
\end{proof}

\begin{proof}[Proof of Theorem \ref{thm: finite upper right Dini derivative for certain periodic tree}]
	The existence of the right derivative of $\theta(p)$ follows from Lemma \ref{lem: remain concave or convex near p_c}.
	The positiveness and finiteness of the right derivatives follow from Proposition \ref{prop: positive lower right Dini derivative} and Lemma \ref{lem: finite upper right Dini derivative} respectively. 
\end{proof}

\section{Concluding remarks and questions}\label{sec: remarks and questions}

\subsection{Remark on Bernoulli site percolation}

For $p\in[0,1]$, if we instead keep each vertex with probability $p$ and remove it otherwise. Call the vertices kept \notion{open vertices} and those removed \notion{closed vertices}. \notion{Bernoulli$(p)$ site percolation} studies the random subgraph $\xi$ of $G$ induced by the open vertices. For Bernoulli site percolation, an edge is call open if and only if its two endpoints are open. When talking about  Bernoulli site percolation, we will use $\mathbb{P}_p^{\mathrm{site}}$ and $\mathbb{P}_p^{\mathrm{site}}$ to stress that. 

\begin{remark}\label{rem: site percolation}
	If the connected graph $G$ has bounded degree, say, with a upper bound $D$, then for  Bernoulli site percolation, the following analogue of Lemma \ref{lem: differential inequality} holds:
	\be\label{eq: differential ineq for site percolation}
	\frac{d }{dp}\mathbb{P}_p^{\mathrm{site}}(x\longleftrightarrow \Lambda_n^c)\geq \frac{1}{1-p}\min(1,\frac{\inf_{S\subset \Lambda_n,x\in S}\varphi_p(x,S)}{D-1})\cdot  [1-\mathbb{P}_p^{\mathrm{site}}(x\longleftrightarrow \Lambda_n^c)]
	\ee
	where $\varphi_p(x,S):=\sum_{y\in S}\sum_{z\notin S, (y,z)\in E}\mathbb{P}_p^{\mathrm{site}}[x\stackrel{S}{\longleftrightarrow} y]$. If one defines $p_{\textnormal{cut,site}}'$ accordingly for site percolation, one can prove $p_{c,\textnormal{site}}=p_{\textnormal{cut,site}}'$ similarly as the bond percolation case. 
\end{remark}
This leads us to the following conjecture:
\begin{conj}\label{conj: site percolation}
	The answers for Question \ref{ques: Kahn} and \ref{ques: Kahn's question for edge cutset} are also positive for Bernoulli site percolation.
\end{conj}

\subsection{Is there an example with $p_{\mathrm{cut,E}}< p_{\mathrm{cut,V}}$?}

We have seen examples with $p_{\mathrm{T,E}}<p_{\mathrm{T,V}}$ (Example \ref{example: p_T,E less than p_T,V}). One can ask the same question for $p_{\mathrm{cut,E}}$ and $p_{\mathrm{cut,V}}$:
\begin{question}
	Is there a locally finite, connected, infinite graph $G$ such that $p_{\mathrm{cut,E}}< p_{\mathrm{cut,V}}$?
\end{question}

In view of  Lemma \ref{lem: p_cut,e less than p_cut,v}, if there is a graph $G$ with $p_{\mathrm{cut,E}}< p_{\mathrm{cut,V}}$, then it must have unbounded degree and for a vertex cutset $\Pi_V$, for ``most" $v\in \Pi_V$ there should be a lot of edges in the corresponding edge cutset $\Pi_E=\Delta S(\Pi_V)$  incident to $v$. One might first want to consider certain 1-dimensional multigraphs. However there is no simple 1-dimensional example; see Proposition \ref{prop: no 1D example for p_cut,e less than p_cut,v}.

\begin{definition}
	Let $(a_n)_{n\geq0}$ be a sequence of positive integers. Let $G=G((a_n)_{n\geq0})$ be the graph with vertex set $V=\mathbb{N}=\{0,1,2,\cdots\}$ and edge set $E=\bigcup_n E_n$, where $E_n=\{e_{n,j}\colon j=1,\ldots,a_n\}$ is the set of $a_n$ parallel edges from $n$ to $n+1$. 
\end{definition}

\begin{proposition}\label{prop: no 1D example for p_cut,e less than p_cut,v}
	There is no sequence 	of positive integers $(a_n)_{n\geq0}$ such that $G=G((a_n)_{n\geq0})$ has the property of  $p_{\mathrm{cut,E}}(G)< p_{\mathrm{cut,V}}(G)$. 
\end{proposition}
We need a technical lemma. 
\begin{lemma}\label{lem: no 1D example}
	There is no  sequence of positive integers $(a_n)_{\geq0}$ that satisfies both
	\be\label{eq: p_2 < p_c}
	\sum_{i=0}^{\infty}(1-p_2)^{a_i}=\infty
	\ee
	and
	\be\label{eq: p_1> p_cut,E}
	a_n\geq \frac{c}{\prod_{i=0}^{n-1}\big[ 1-(1-p_1)^{a_i} \big]}, \forall\, n\geq1.
	\ee
	for some constants $0<p_1<p_2<1$ and  $c>0$. 
\end{lemma}
\begin{proof}[Proof of Proposition \ref{prop: no 1D example for p_cut,e less than p_cut,v} assuming Lemma \ref{lem: no 1D example}]
	Notice that 
	\be\label{eq: two point function}
	\mathbb{P}_p[0\longleftrightarrow n]=\prod_{i=0}^{n-1}\big[ 1-(1-p)^{a_i} \big]
	\ee
	
	Thus
	\be\label{eq: p < p_c}
	p<p_{\mathrm{c}}  \,\,\Rightarrow \,\, \lim_{n\to\infty}\prod_{i=0}^{n-1}\big[ 1-(1-p)^{a_i} \big]=0\,\, \Leftrightarrow \,\, \sum_{i=0}^{\infty}(1-p)^{a_i}=\infty
	\ee
	and 
	\be
	\sum_{i=0}^{\infty}(1-p)^{a_i}<\infty  \,\,\Leftrightarrow \,\,\lim_{n\to\infty}\prod_{i=0}^{n-1}\big[ 1-(1-p)^{a_i} \big]>0  \,\,\Rightarrow \,\, p\geq p_{\mathrm{c}}
	\ee
	
	Notice that the minimal vertex cutsets are $\{\Pi_n\}$ where $\Pi_n=n$ and the minimal edge cutsets are $E_n=\{e_{n,j}\colon j=1,\cdots,a_n\}$, the $a_n$ parallel edges from $n$ to $n+1$.

	Hence $\mathbb{E}_p[|C(0)\cap \Pi_n|]=\mathbb{P}_p[0\longleftrightarrow n]=\prod_{i=0}^{n-1}\big[ 1-(1-p)^{a_i} \big]$. Thus $p_{\mathrm{cut,V}}=p_{\mathrm{c}}$. 
	
	Suppose there is some sequence $(a_n)_{n\geq0}$ such that $p_{\mathrm{cut,E}}<p_{\mathrm{cut,V}}$. Pick $p_1,p_2$ such that $p_{\mathrm{cut,E}} <p_1<p_2< p_{\mathrm{cut,V}}$. 
	
	Since $p_2< p_{\mathrm{cut,V}}=p_{\mathrm{c}}$, by \eqref{eq: p < p_c} one has that 
	\[
	\sum_{i=0}^{\infty}(1-p_2)^{a_i}=\infty
	\]
	
	Since we choose $p_1>p_{\mathrm{cut,E}}$ and noting that $\mathbb{E}_p[|C(0)\cap E_n|]=pa_n\mathbb{P}_p[0\longleftrightarrow n]$, one has that
	\be
	\inf_{n}  p_1a_n\prod_{i=0}^{n-1}\big[ 1-(1-p_1)^{a_i} \big]>0
	\ee
	i.e., there exists $c>0$ s.t.\ 
	\[
	a_n\geq \frac{c}{\prod_{i=0}^{n-1}\big[ 1-(1-p_1)^{a_i} \big]}, \forall\, n\geq1.
	\]
	But this contradicts with Lemma \ref{lem: no 1D example} and hence Proposition \ref{prop: no 1D example for p_cut,e less than p_cut,v} holds. 
\end{proof}

\begin{proof}[Proof of Lemma \ref{lem: no 1D example}]
	First we reduce to the case of increasing sequence. 
	If there is some sequence  satisfies both \eqref{eq: p_2 < p_c} and \eqref{eq: p_1> p_cut,E} for some $0<p_1<p_2<1$,  
	then 
	\[
	\prod_{i=0}^{n-1}\big[ 1-(1-p_1)^{a_i} \big]\leq \prod_{i=0}^{n-1}\big[ 1-(1-p_2)^{a_i} \big] \stackrel{\eqref{eq: p_2 < p_c}}{\to} 0\textnormal{ as }n\to \infty.
	\]	
	Thus by \eqref{eq: p_1> p_cut,E}, one has that $a_n\to \infty$. In the following for simplicity we write $p=p_1$.	
	
	Now we consider the sequence $(a_n')$, the rearrangement of $a_n$ in the non-decreasing order.  Obviously, $(a_n')$ also satisfies \eqref{eq: p_2 < p_c}. As for \eqref{eq: p_1> p_cut,E}, 
	let $m=m(n)$ be the last index such that $a_m\leq a_n'$, i.e., $m=\max\{k\colon a_k\leq a_n'\}$. Obviously $m\geq n$.  Since $a_k\to\infty$ as $k\to\infty$, $m<\infty$. 
	Then we claim that there exists a constant $c_0>0$ such that
	\be\label{eq: key inequality for increasing reduction}
	a_n'\prod_{i=0}^{n-1}\big[ 1-(1-p)^{a_i'} \big]
	\geq c_0a_m\prod_{i=0}^{m-1}\big[ 1-(1-p)^{a_i} \big] \stackrel{\eqref{eq: p_1> p_cut,E}}{\geq} cc_0,
	\ee
	%		where the second inequality $\stackrel{?}{\geq}$ is due to the fact that the multi-set $\{a_0',\cdots,a_{n-1}'\}$ must be contained in $\{a_0,\cdots,a_{m-1}\}$  by our choice of $m$: more specifically\\ 
	
	Write $A=\{v_i\colon v_1<v_2<\cdots\}$  for the all the values of the sequence $(a_n)$. For each $v\in A$, let $N(v)=|\{ j\colon a_j=v \}|\geq1$ be the number of times  the sequence taking the value $v$. 	
	
	\textbf{Case one:} $a_m=a_n'$. By the definition of  $a_n'$, we assume $a_n'=v_k$ and then 
	\be\label{eq: 0.8}
	\prod_{i=0}^{n-1}\big[ 1-(1-p)^{a_i'} \big]\geq \big[ 1-(1-p)^{v_k} \big]^{N(v_k)-1}\times \prod _{i=1}^{k-1}\big[ 1-(1-p)^{v_i} \big]^{N(v_i)}
	\ee
	By the choice of $m$, the multi-set $\{a_0,\cdots,a_{m-1}\}$ contains at least $N(v_k)-1$'s $v_k$ and all the $N(v_i)$'s $v_i$ for $i<k$. Hence 
	\be\label{eq: 0.9}
	\big[ 1-(1-p)^{v_k} \big]^{-1}\times \prod _{i=1}^{k}\big[ 1-(1-p)^{v_i} \big]^{N(v_i)}\geq \prod_{i=0}^{m-1}\big[ 1-(1-p)^{a_i} \big].
	\ee
	By $a_n'= a_m$ and the above two inequalities \eqref{eq: 0.8},\eqref{eq: 0.9} we have \eqref{eq: key inequality for increasing reduction} for all $c_0\leq1$. 
	
	\textbf{Case two:} $a_m<a_n'$, say $a_n'=v_k$ and $a_m=v_j$ for some $j<k$.
	Then 	by the choice of $m$, the multi-set $\{a_0,\cdots,a_{m-1}\}$ contains at least $N(v_j)-1$'s $v_j$ and all the  other $N(v_i)$'s $v_i$ for $i\leq k, i\neq j$. Hence 
	\be\label{eq: 0.10}
	\big[ 1-(1-p)^{v_j} \big]^{-1}\times \prod _{i=1}^{k}\big[ 1-(1-p)^{v_i} \big]^{N(v_i)}\geq \prod_{i=0}^{m-1}\big[ 1-(1-p)^{a_i} \big].
	\ee
	
	By \eqref{eq: 0.8} and \eqref{eq: 0.10} we have that 
	\[
	\frac{ 	a_n'\prod_{i=0}^{n-1}\big[ 1-(1-p)^{a_i'} \big] }{ a_m\prod_{i=0}^{m-1}\big[ 1-(1-p)^{a_i} \big] }  \geq \frac{v_k \big[ 1-(1-p)^{v_k} \big]^{-1} }{v_j \big[ 1-(1-p)^{v_j} \big]^{-1}} \geq p\frac{v_k}{v_j}\geq p.
	\]
	where in the second inequality we use the fact that 
	\[
	f(x)=\frac{x}{1-(1-p)^x}\in (x,\frac{x}{p}],x\geq1.
	\]
	Hence in this case \eqref{eq: key inequality for increasing reduction} holds with $c_0=p=p_1$.
	
	Combining the two cases one has that \eqref{eq: key inequality for increasing reduction} holds with $c_0=p=p_1$.

	By the reduction, we can assume $(a_n)_{n\geq0}$ is increasing. Thus there is a strictly increasing sequence  $(n_j)$ such that for $n\in [n_j,n_{j+1}-1]$, $a_n=v_j$. In particular, 
	\eqref{eq: p_2 < p_c} becomes 
	\be\label{eq: p_2 < p_c for increasing case}
	\sum_{j=1}^{\infty}(n_{j+1}-n_j)(1-p_2)^{v_j}=\infty 
	\ee	
	and 
	\eqref{eq: p_1> p_cut,E} becomes (only needs to look at times $n_{j+1}-1$)
	\be\label{eq: p_1>p_cut,E for increasing case}
	v_{j}\geq \frac{c\big[ 1-(1-p_1)^{v_j}\big]}{\prod_{i=1}^{j}\big[ 1-(1-p_1)^{v_i} \big]^{n_{i+1}-n_i}}
	\ee	
	
	By \eqref{eq: p_1>p_cut,E for increasing case}  one has that 
	\[
	\frac{c}{v_j\big[ 1-(1-p_1)^{v_j} \big]^{n_{j+1}-n_j-1}}\leq \prod_{i=1}^{j-1}\big[ 1-(1-p_1)^{v_i} \big]^{n_{i+1}-n_i} \leq 1.
	\]	
	Hence 
	\[
	v_j\big[ 1-(1-p_1)^{v_j} \big]^{n_{j+1}-n_j-1} \geq c.
	\]	
	Taking logarithm one has that 
	\[
	\log v_j+(n_{j+1}-n_j-1)\log [1-(1-p_1)^{v_j}] \geq\log c
	\]
	Hence 
	\[
	n_{j+1}-n_j-1\leq \frac{\log c-\log v_j}{\log [1-(1-p_1)^{v_j}]}\leq \frac{c'\log v_j}{(1-p_1)^{v_j}},\textnormal{ when }v_j>1.
	\]
	But this contradicts with \eqref{eq: p_2 < p_c for increasing case}: (noting $\{v_j\}$ is a strictly increasing subsequence of $\mathbb{N}$ and $1-p_2<1-p_1$)
	\[
	\sum_{j=1}^{\infty}(n_{j+1}-n_j)(1-p_2)^{v_j}\leq (n_2-n_1)(1-p_2)^{v_1}+\sum_{j=2}^{\infty} (1-p_2)^{v_j}+\sum_{j=2}^{\infty} \frac{c'\log v_j}{(1-p_1)^{v_j}} (1-p_2)^{v_j}<\infty.
	\]
	
	This contradiction implies Lemma \ref{lem: no 1D example}.
\end{proof}

\section*{Acknowledgments}
We thank Russ Lyons for many helpful discussions and support. Most of the work was done while the author was a graduate student at Indiana University.  We also thank Wai-Kit Yeung for discussions and the reference on semi-algebraic functions. We are also grateful to the anonymous referee for his/her careful reading and helpful comments which improve the presentation a lot.

%% General advice from Russell Lyons. %%
%% After solving a problem, one could look at several problems. %%
%% First, could the proof be simplified? %%
%% Second, could the result refine other related results? (application and relation to other results) %%
%% Third, could one raise new  questions?%%

\bibliography{pcpcut}
\bibliographystyle{plain}

\end{document}